\newcommand{\dsp}{\renewcommand{\baselinestretch}{1.0}}

\newcommand{\lr}{\longrightarrow}
\newcommand{\cd}{\cdots}

\newcommand{\N}{\mathbb{N}}
\newcommand{\C}{\mathbb{C}}
\newcommand{\R}{\mathbb{R}}

\newcommand{\Z}{\mathbb{Z}}
\newcommand{\T}{\mathbb{T}}
\newcommand{\qq}{\quad}
\newcommand{\af}{\alpha}
\newcommand{\one}{\textbf{1}}

\documentclass[pdflatex,fleqn]{article}
\usepackage{}
\usepackage{xcolor}
\usepackage{mathrsfs}
\usepackage{amsfonts}
\usepackage{amssymb}
\usepackage{latexsym}
\usepackage{amsmath}
\usepackage{amscd}
\usepackage{graphicx}
\usepackage{hyperref}
\usepackage{hypernat}
\usepackage{cite}
\usepackage{mathtools}

\usepackage{geometry,amsthm,graphics,amssymb,amsmath,enumerate,latexsym,tabularx,shapepar}
\usepackage[all,2cell,pdftex]{xy}\UseAllTwocells\SilentMatrices

\usepackage{extarrows}
\date{}

\dsp

\begin{document}


\centerline{\large \bf Hausdorffifized  algebraic $K_1$ group and invariants for $C^*$-algebras with the ideal property }

\vspace{3mm}

\centerline{Guihua Gong, Chunlan Jiang, and Liangqing Li}

\vspace{3mm}

\centerline{Dedicated to the  memory of Professor Ronald G. Douglas}







\noindent{\bf Abstract}
 A $C^*$-algebra $A$ is said to have the ideal property if each closed two-sided ideal of $A$ is
generated by the projections inside the ideal, as a closed two sided ideal. $C^*$-algebras with the ideal property are   generalization and  unification of real rank zero $C^*$-algebras and unital simple $C^*$-algebras. It is long to be expected that an invariant (see [Stev] and [Ji-Jiang], [Jiang-Wang] and [Jiang1]) , we call it $Inv^0(A)$ (see the introduction),  consisting of scaled ordered total
$K$-group $(\underline{K}(A), \underline{K}(A)^{+},\Sigma A)_{\Lambda}$ (used in the real rank zero case), the tracial
state space $T(pAp)$ of cutting down algebra $pAp$ as part of Elliott invariant of $pAp$ (for each $[p]\in\Sigma A$) with a certain compatibility, is the complete invariant for certain well behaved class of $C^*$-algebras with the ideal property (e.g., $AH$ algebras with no dimension growth). In this paper, we will  construct two non isomorphic  $A\mathbb{T}$ algebras $A$ and $B$ with the ideal  property such that $Inv^0(A)\cong Inv^0(B)$. The invariant to differentiate the two algebras is the Hausdorffifized  algebraic $K_1$-groups $U(pAp)/\overline{DU(pAp)}$ (for each $[p]\in\Sigma A$) with a certain compatibility condition. It will be proved in [GJL] that, adding this new ingredients, the invariant will become the complete invariant for $AH$ algebras (of no dimension growth) with the ideal property.







\noindent\textbf{\S1. Introduction}

A $C^*$-algebra $A$ is called an $AH$ algebra (see {Bl])  if it is the inductive limit $C^{*}$-algebra
of
$$A_{1}\xrightarrow{\phi_{1,2}}A_{2}\xrightarrow{\phi_{2,3}}A_{3}\longrightarrow\cdots\longrightarrow A_{n}\longrightarrow\cdots$$ with
 $A=\lim\limits_{n\to \infty}(A_{n}=\bigoplus\limits_{i=1}\limits^{t_{n}}P_{n,i}M_{[n,i]}(C(X_{n,i}))P_{n,i}, \phi_{n,m})$, where $X_{n,i}$ are compact metric spaces, $t_{n}$ and $[n,i]$ are
positive integers, and $P_{n,i}\in M_{[n,i]}(C(X_{n,i}))$ are projections. An $AH$ algebra is called of no
dimension growth, if one can choose the spaces $X_{n,i}$ such that $sup_{n,i}dim(X_{n,i})<+\infty$. If all the spaces $X_{n,i}$ can be chosen to be the single point space $\{pt\}$, then $A$ is called an $AF$ algebra. If all the spaces can be chosen to be the interval $[0,1]$ (or circle $\T=\{z\in \C: |z|=1\}$, respectively) , then $A$ is called an $AI$ algebra (or $A\T$ algebras, respectively).

In 1989, G. Elliott (see [Ell1]) initiated  the classification program by classying all real rank zero $A\T$ algebras (without the condition of simplicity) and he conjectured that the scaled ordered $K_*$ group \\ $(K_*(A), K_*(A)^+, \Sigma A)$ , where $K_*(A)=K_0(A)\oplus K_1(A)$, is a complete invariant for separable nuclear $C^*$-algebras of real rank zero and stable rank one. In 1993,  Elliott (see [Ell2]) successfully classified all unital simple $AI$ algebras by the so called Elliott invariant $Ell(A)=(K_0(A), K_0(A)^+, \Sigma A, K_1(A), TA,\rho_A)$, where $TA$ is the space of
all unital traces on $A$, and $\rho_A$ is the nature map from $K_0(A)$ to $Aff TA$ (the ordered Banach space of all affine maps from $TA$ to $\R$).

In 1994, the first named author (see [G1])  constructed two non isomorphic (not simple) real rank zero $AH$ algebras (with 2-dimensional local spectra) $A$ and $B$ such that $(K_*(A), K_*(A)^+, \Sigma A)\cong (K_*(B), K_*(B)^+, \Sigma B)$, which disproved the conjecture of Elliott for $C^*$-algebras of real rank zero and stable rank one. This result lead to a sequence of research by Dadarlat-Loring, Eilers (see [DL1-2]. [Ei]) end up with Dadarlat-Gong's  complete classification  (see [DG]) of real rank zero $AH$ algebras by scaled ordered total K-theory $(\underline{K}(A), \underline{K}(A)^{+},\Sigma A)_{\Lambda}$, where $\underline{K}(A)=K_*(A)\oplus \bigoplus_{p=2}^{\infty} K_*(A, \Z/p\Z)$ and $\Lambda$ is the system of Bockstein operations (also see [D1-2],  [EG1-2],[EGLP], [EGS], [G1-4], [GL] and [Lin1-3]). In [EGL1], Elliott-Gong-Li completely classified simple $AH$ algebras of no dimension growth by Elliott invariant (also see [Ell3], [EGL2], [EGJS], [G5], [Li1-5], [Lin4], [NT] and [Thm1-2]).  A natural  generalization and unification of real rank zero $C^*$-algebras and unital simple $C^*$-algebras is the class of $C^*$-algebras with the ideal property:  each closed two-sided ideal is
generated by the projections inside the ideal, as a closed two sided ideal. It is long to be expected that a combination of scaled ordered total K-theory (used in the classification of real rank zero $C^*$-algebras) and the Elliott invariant (used in the the classification of simple $C^*$-algebras), including tracial state spaces $T(pAp)$---part of Elliott invariant  of cutting down algebras $\{pAp\}_{[p]\in \Sigma A}$ with comptibility conditions, called $Inv^0(A)$ (see 2.18 of [Jiang1]),  is a
 complete invariant for certain well behaved (e.g., ${\cal Z}$-stable, where ${\cal Z}$ is the Jiang-Su algebra of [JS]) $C^*$-algebras with the ideal property (see [Stev], [Pa], [Ji-Jiang],[Jiang-Wang], [Jiang1]).

The main purpose of this paper is to construct two unital ${\cal Z}$-stable  A$\T$ algebras $A$ and $B$ with the ideal property such that $Inv^0(A)\cong Inv^0(B)$, but $A\not\cong B$.
 The invariant to distinguish  these two $C^*$-algebras is the Hausdorffifized  algebraic $K_1$ groups $U(pAp)/\overline{DU(pAp)}$ of the cutting down algebra  $pAp$ (for each element $x\in \Sigma A$, we chose one projection $p\in A$ such that $[p]=x$) with a certain compatibility condition, where $DU(A)$ is the group generated  by commutators $\{uvu^*v^*~|~u, v\in U(A)\}$.
 In this paper, we will  introduce the invariant $Inv'(A)$ and its simplified  version $Inv(A)$, by adding these new ingredients--the Housdorffized algebraic $K_1$ groups of cutting down algebras with compatibility conditions,  to $Inv^0(A)$.

 In [GJL], we will prove that $Inv(A)$ is a complete invariant for $AH$ algebras (of no dimension growth) with the ideal property.

 Let us point out that for the above  $C^*$-algebras $A$ and $B$, we have that $Cu(A)\cong Cu(B)$ and $Cu(A\otimes C(S^1))\cong Cu(B\otimes C(S^1))$. That is, the new invariant can not  be detected by Cuntz semigroup.

 In section 2, we will define $Inv(A)$ and discuss its properties. These properties will be used in [GJL].  In section 3, we will present the construction of $A\T$ algebras $A$ and $B$ with the ideal property such that $Inv(A)\not\cong Inv(B)$, but $Inv^0(A)\cong Inv^0(B)$.








\noindent\textbf{2. The invariant}

In this section, we will recall the definition of $Inv^0(A)$ from [Jiang1] (also see [Stev],  [Ji-Jiang], [Jiang-Wang]), and then introduce the invariant $Inv(A)$. Furthermore, we will discuss the properties of $Inv(A)$ in the context of $AH$ algebras and $A\mathcal{HD}$ algebras (for definition of $A\mathcal{HD}$ algebras, see 2.3 below), which are used in [GJL].

\noindent\textbf{2.1.}~~In the notation for an inductive limit system $\lim (A_{n},\phi_{n,m})$, we understand that $$\phi_{n,m}=\phi_{m-1,m}\circ\phi_{m-2,m-1}\circ\cdots\circ\phi_{n,n+1},$$
where all $\phi_{n,m}:A_{n}\rightarrow A_{m}$ are homomorphisms.

We shall assume that, for any summand $A^{i} _{n}$ in the direct sum $A_{n}=\bigoplus^{t_{n}}_{i=1}A^{i} _{n}$, necessarily,
$\phi_{n,n+1}(\textbf{1}_{A^{i} _{n}})\neq0$, since, otherwise, we could simply delete $A^{i} _{n}$ from $A_{n}$, without changing the limit algebra.

If $A_{n}=\bigoplus_{i}A^{i} _{n}$, $A_{m}=\bigoplus_{j}A^{j} _{m}$, we use $\phi^{i,j}_{n,m}$ to denote the partial map of $\phi_{n,m}$ from the
$i$-th block $A^{i} _{n}$ of $A_{n}$ to the $j$-th block $A^{j} _{m}$ of $A_{m}$. Also, we use $\phi_{n,m}^{-,j}$ to denote the partial
map of $\phi_{n,m}$ from $A_{n}$ to $A_{m}^{j}$. That is,  $\phi_{n,m}^{-,j}=\bigoplus\limits_{i}\phi_{n,m}^{i,j}=\pi_{j}\phi_{n,m}$,
where $\pi_{j}:A_{m}\rightarrow A_{m}^{j}$ is the canonical projection. Some times, we also use $\phi_{n,m}^{i,-}$ to denote $\phi_{n,m}|_{A_n^i}:~ A_n^i \to A_m$.

\noindent\textbf{2.2.}~~As in [EG2], let $T_{\uppercase\expandafter{\romannumeral2},k}$ be the $2$-dimensional connected simplicial complex with $H^{1}(T_{\uppercase\expandafter{\romannumeral2},k})=0$ and $H^{2}(T_{\uppercase\expandafter{\romannumeral2},k})=\mathbb{Z}/k\mathbb{Z}$, and
let $I_{k}$ be the subalgebra of $M_{k}(C[0,1])=C([0,1],M_{k}(\mathbb{C}))$ consisting of all functions $f$ with the properties $f(0)\in \mathbb{C}\cdot \textbf{1}_{k}$ and $f(1)\in \mathbb{C}\cdot \textbf{1}_{k}$
(this algebra is called an Elliott dimension drop interval algebra). Denoted by  $\mathcal{HD}$ the class of algebras consisting of direct sums of the building
blocks of the forms $M_{l}(I_{k})$ and $PM_{n}(C(X))P$, with $X$ being one of the spaces $\{pt\}$, $[0,1]$, $S^{1}$, and $T_{\uppercase\expandafter{\romannumeral2},k}$, and with $P\in M _{n}(C(X))$ being a
projection. (In [DG], this class is denoted by $SH(2)$, and in [Jiang1], this class is denoted by $\mathcal{B}$). We will call a $C^{*}$-algebra an $A\mathcal{HD}$
algebra, if it is an inductive limit of the algebras in $\mathcal{HD}$.

For each basic building block $A=PM_{n}(C(X))P$, where $X=\{pt\}, [0,1], S^1, T_{\uppercase\expandafter{\romannumeral2},k}$, or $A=M_l(I_k)$, we have $K_0(A)= \Z $ or $\Z/k\Z$ (for the case $A=P(M_n(C(T_{\uppercase\expandafter{\romannumeral2},k}))P$). Hence there is a natural map $rank: K_0(A) \to \Z$. This map also gives a map from $\{p\in (M_{\infty}(A)):~p ~\mbox{is a projection} \}$ to $\Z_+$. For example, if $p\in A=PM_{n}(C(X))P$, then $rank (p)$ is the rank of projection $p(x)\in P(x)M_n(\C)P(x)\cong M_{rank (P)}(\C)$ for any $x\in X$; and if $p\in A=M_l(I_k)$, then $rank (p)$ is the rank of projection $p(0)\in M_l(\C)$. (Note that we regard $p(0)$ in $M_l(\C)\cong \one_k\otimes M_l(\C)$ (not regard it  in $ M_{lk}(\C)$).)


\noindent\textbf{2.3.}~~By $A\mathcal{HD}$ algebra, we mean the inductive limit of
$$A_{1}\xrightarrow{\phi_{1,2}}A_{2}\xrightarrow{\phi_{2,3}}A_{3}\longrightarrow\cdots\longrightarrow\cdots,$$
where $A_{n}\in \mathcal{HD}$ for each $n$.

For  an $A\mathcal{HD}$ inductive limit $A\!=\!\lim (A_n, \phi_{nm})$, we write $A_n\!=\!\oplus_{i=1}^{t_n}A_n^i$, where
$A_n^i=P_{n,i}M_{[n,i]}(C(X_{n,i}))P_{n,i}$ or $A_n^i=M_{[n,i]}(I_{k_{n,i}})$. For convenience, even for a block $A_n^i=M_{[n,i]}(I_{k_{n,i}})$, we still use  $X_{n,i}$ for $Sp(A_n^i)=[0,1]$---that is, $A_n^i$ is regarded as a homogeneous algebra or a  sub-homogeneous algebra over $X_{n,i}$.

\noindent\textbf{2.4.}~~In [GJLP1-2], joint with Cornel Pasnicu, the authors proved the reduction theorem for $AH$ algebras with the ideal property provided that the inductive limit systems
have no dimension growth. That is, if $A$ is an inductive limit of $A_{n}=\bigoplus A^{i}_{n}=\bigoplus P_{n,i}M_{[n,i]}C(X_{n,i})P_{n,i}$ with $sup_{n,i}dim(X_{n,i})<+\infty$, and if we further assume that $A$ has the ideal property, then $A$ can be rewritten as an  inductive limit of $B_{n}=\bigoplus B^{j}_{n}=\bigoplus Q_{n,j}M_{\{n,j\}}C(Y_{n,i})Q_{n,j}$, with $Y_{n,i}$ being one of $\{pt\}$, $[0,1]$, $S^{1}$, $T_{\uppercase\expandafter{\romannumeral2},k}$, $T_{\uppercase\expandafter{\romannumeral3},k}$, $S^{2}$. In turn, the
 second author proved in [Jiang2] (also see [Li4]), that the above inductive limit can be rewritten as the inductive limit of the direct sums of homogeneous algebras
over $\{pt\}$, $[0,1]$, $S^{1}$, $T_{\uppercase\expandafter{\romannumeral2},k}$ and $M_{l}(I_{k})$. Combining these two results, we know that all $AH$ algebras of no dimension growth with the ideal
property are $A\mathcal{HD}$ algebras. Let us point out that,  as proved in [DG], there are real rank zero $A\mathcal{HD}$ algebras which are not $AH$ algebras.

\noindent\textbf{2.5.}~~Let A be a $C^{*}$-algebra. $K_{0}(A)^{+}\subset K_{0}(A)$ is defined to be the semigroup of $K_{0}(A)$ generated
by $[p]\in K_{0}(A)$, where $p\in M_{\infty}(A)$ are projections. For all $C^{*}$-algebras considered in this paper, for
example, $A\in \mathcal{HD}$, or $A$ is an $A\mathcal{HD}$ algebra, or $A=B\otimes C(T_{\uppercase\expandafter{\romannumeral2},k}\times S^{1})$, where B is an $\mathcal{HD}$ or $A\mathcal{HD}$
algebra, we  always have
$$(*)\qq\qq\qq\qq\qq\qq K_{0}(A)^{+}\bigcap (-K_{0}(A)^{+})=\{0\}~~~\mbox{and}~~~K_{0}(A)^{+}-K_{0}(A)^{+}=K_{0}(A).\qq\qq\qq$$
Therefore $(K_{0}(A),K_{0}(A)^{+})$ is an ordered group. Define $\Sigma A \subset K_{0}(A)^{+}$ to be
$$\Sigma A = \{[p]\in K_{0}(A)^{+}, p~ is ~a ~projection ~in ~A\}.$$
Then $(K_{0}(A),K_{0}(A)^{+},\Sigma A)$ is a scaled ordered group. (Note that for purely infinite $C^*$ algebras or stable projectionless $C^*$algebras, the above condition $(*)$ does not hold.)

\noindent\textbf{2.6.}~~Let $\underline K (A)=K_{*}(A)\bigoplus \big(\bigoplus_{k=2}^{+\infty}K_{*}(A,\mathbb{Z}/ k\mathbb{Z})\big)$ be as in [DG].
Let $\wedge$ be the Bockstein operation on $\underline K (A)$(see 4.1 of [DG]). It is well known that
$K_{*}(A,Z\oplus \mathbb{Z}/ k\mathbb{Z})=K_{0}(A\otimes C(W_{k}\times S^{1})),$
where $W_{k}=T_{II,k}$.

As in [DG], let
$K_{*}(A,Z\oplus \mathbb{Z}/ k\mathbb{Z})^{+}=K_{0}(A\otimes C(W_{k}\times S^{1})^{+})$
and let $\underline K (A)^{+}$ be the semigroup generated by $\{K_{*}(A,\mathbb{Z}\oplus\mathbb{Z}/ k\mathbb{Z})^{+},k=2,3,\cdots\}$.

\noindent\textbf{2.7.}~~Let $Hom_{\wedge} (\underline K (A),\underline K (B))$ be the set of homomorphisms between $\underline K (A)$ and $\underline K (B)$
compatible with the Bockstein operations $\wedge$. There is   a surjective map (see [DG])
$$\Gamma: KK(A,B)\rightarrow Hom_{\wedge}(\underline K (A),\underline K (B)).$$
Following  R\o rdam(see [R]), we denote $KL(A,B)\triangleq KK(A,B)/Pext(K_*(A), K_{*+1}(B))$, where \\$Pext(K_*(A), K_{*+1}(B))$  is identified with  $ker\;\Gamma$ by [DL2].  The triple
 $(\underline{K}(A),\underline{K}(A)^{+},\Sigma A)$ is part of our invariant.
For two $C^{*}$-algebras $A$ and $B$, by a ``homomorphism"
$$\alpha:~(\underline{K}(A),\underline{K}(A)^{+},\Sigma A) \to (\underline{K}(B),\underline{K}(B)^{+},\Sigma B),$$
 we mean a system of maps:
$$\alpha^{i}_{k}: ~K_{i}(A,\mathbb{Z}/k\mathbb{Z})\longrightarrow K_{i}(B,\mathbb{Z}/k\mathbb{Z}),~~i=0,1,~~k=0,2,3,\cdots$$ which are compatible with the Bockstein operations  and $\alpha=\oplus_{k,i}\alpha^{i}_{k}$ satisfies $\alpha(\underline{K}(A)^{+})\subset\underline{K}(B)^{+}$. And finally,  $\alpha^{0}_{0}(\Sigma A)\subset\Sigma B$.

\noindent\textbf{2.8.}~~ For a unital $C^{*}$-algebra $A$, let $TA$ denote the space of tracial states of $A$, i.e., $\tau\in TA$ if and only if $\tau$ is a positive linear map from $A$ to $\mathbb{C}$ with $\tau(xy)=\tau(yx)$, and $\tau(\textbf{1})=1$. Endow $TA$ with the weak-* topology, that is, for any net $\{\tau_{\af}\}_{\af} \subset TA$ and $\tau\in TA$, $\tau_{\af}\to \tau$ if and only if $\lim_{\af}\tau_{\af}(x)=\tau (x)$ for any $x\in A$. Then $TA$ is a compact Hausdorff space with convex structure, that is, if $\lambda\in [0,1]$ and $\tau_1, \tau_2\in TA$, then $\lambda \tau_1+(1-\lambda)\tau_2\in TA$.  $AffTA$ is the collection of all continuous  affine maps from $TA$ to $\mathbb{R}$, which is a real Banach space with $\|f\|=\mbox{sup}_{\tau \in TA}|f(\tau)|$.
Let $(AffTA)_{+}$  be the subset of $AffTA$ consisting of all nonnegative affine functions.
An element
$\textbf{1}\in AffTA$, defined by $\textbf{1}(\tau)=1$ for all $\tau\in TA$,  is called the order unit (or scale)  of $AffTA$. Note that any $f\in AffTA$ can be written as $f=f_+-f_-$ with $f_1, f_2 \in AffTA_+$, $\|f_1\|\leq \|f\|$  and $\|f_2\|\leq \|f\|$. Therefore  ($AffTA$, $(AffTA)_{+}$, \textbf{1}) forms a scaled ordered real Banach space. If $\phi: AffTA \to AffTB$ is a unital positive linear map, then $\phi$ is bounded and therefore continuous.

There is a natural homomorphism $\rho_A: K_0(A) \to AffTA$ defined by $\rho_A([p])(\tau)=\sum_{i=1}^n \tau(p_{ii})$ for $\tau\in TA$ and $[p]\in K_0(A)$ represented by projection $p=(p_{ij})\in M_n(A)$.

Any unital homomorphism $\phi: A\longrightarrow B$ induces a continuous affine map $T\phi: TB\longrightarrow TA$, which, in turn,  induces a unital  positive linear map $AffT\phi: AffTA\longrightarrow AffTB.$

If $\phi: A\longrightarrow B$ is not unital, we still use $AffT\phi$ to denote the unital  positive linear map $$AffT\phi: AffTA\longrightarrow AffT(\phi(\textbf{1}_A)B\phi(\textbf{1}_A))$$
by regarding $\phi$ as the unital homomorphism from $A$ to $\phi(\textbf{1}_A)B\phi(\textbf{1}_A)$---that is, for any $l\in AffTA$ represented by $x\in A_{s.a}$ as $l(t)=t(x)$ for any $t\in TA$, we define
$$\big((AffT\phi)(l)\big)(\tau)=\tau(\phi(x))~~\mbox{for any} ~~\tau \in T(\phi(\one_A)B\phi(\one_A)),$$ where $\phi(x)$ is regarded as an element in $\phi(\one_A)B\phi(\one_A)$. In the above equation, if we regard $\phi(x)$ as element in $B$ (rather than in $\phi(\one_A)B\phi(\one_A)$), the homomorphism $\phi$ also induces a positive linear map, denoted by $\phi_T$ to avoid the confusion,  from $AffTA$ to $AffTB$---that is
for the $l$ as above, $$\big((\phi_T)(l)\big)(\tau)=\tau(\phi(x))~~\mbox{for any} ~~\tau \in T(B),$$ where $\phi(x)$ is now regarded as an element in $B$.  But this map  will not preserve the
unit $\textbf{1}$. It has the property that $\phi_T(\textbf{1}_{AffTA})\leq\textbf{1}_{AffTB}$.

In this paper,  we will often use the notation $\phi_T$ for the following situation: If $p_1<p_2$ are two projections in $A$, and $\phi=\imath:~p_1Ap_1 \longrightarrow p_2Ap_2 $ is the inclusion, then $\imath_T$ will denote the (not necessarily unital) map from $AffT(p_1Ap_1)$ to $AffT(p_2Ap_2)$ induced by $\imath$.

\noindent\textbf{2.9.}~~If $\alpha: (\underline{K}(A),\underline{K}(A)^{+},\Sigma A)\longrightarrow (\underline{K}(B),\underline{K}(B)^{+},\Sigma B)$ is a homomorphism as in 2.7, then for each projection $p\in A$, there is a projection $q\in B$ such that $\alpha([p])=[q].$

Since $I_k$ has stable rank one and the spaces $X$ involved in the definition of $\mathcal{HD}$ class (see $PM_n(C(X))P$ in 2.2) are of dimension at most two, we know that  for all $C^{*}$-algebras $A$ considered in this paper---$\mathcal{HD}$ class or $A\mathcal{HD}$ algebra,  the following statement is true: If $[p_{1}]=[p_{2}]\in K_{0}(A)$, then there is a unitary $u\in A$ such that  $up_{1}u^{*}=p_{2}$. Therefore,  both $AffT(pAp)$ and $AffT(qBq)$ depend only on the
classes $[p]\in K_{0}(A)$ and $[q]\in K_{0}(B)$, respectively. Furthermore, if $[p_{1}]=[p_{2}]$, then the identification of $AffT(p_{1}Ap_{1})$ and $AffT(p_{2}Ap_{2})$ via the
unitary equivalence $up_{1}u^{*}=p_{2}$ is canonical---that is, it does not depend on the choice of unitary $u$. For classes $[p]\in\Sigma A(\subset K_{0}(A)^{+}\subset K_{0}(A))$, we will also take $AffT(pAp)$ as  part of our invariant. We will consider a system of unital positive linear maps $$\xi^{p,q}: AffT(pAp)\longrightarrow AffT(qBq)$$
associated with all pairs of  two classes $[p]\in\Sigma A$ and $[q]\in\Sigma B$, with $\alpha([p])=[q]$.
Such system of maps is said to be compatible if for any $p_{1}\leq p_{2}$ with $\alpha([p_{1}])=[q_{1}]$, $\alpha([p_{2}])=[q_{2}]$, and $q_{1}\leq q_{2}$, the following diagram commutes $$\CD
  AffT(p_{1}Ap_{1}) @>\xi^{p_{1},q_{1}}>> AffT(q_{1}Bq_{1}) \\
  @V\imath_T  VV @V \imath_T VV  \\
  AffT(p_{2}Ap_{2}) @>\xi^{p_{2},q_{2}}>> AffT(q_{2}Bq_{2}),
\endCD\;\;\;\;\;\;\;\;\;~~~~~~~~~~~~~~~~~~~~~(2.A)$$
where the verticle maps are induced by the inclusions.
(See [Ji-Jiang] and [Stev].)

\noindent\textbf{2.10.}~~In this paper, we will denote $$(\underline{K}(A), \underline{K}(A)^{+}, \Sigma A,~ \{AffT(pAp)\}_{[p]\in\Sigma A})$$ by $Inv^{0}(A)$, where $AffT(pAp)$ are scaled ordered Banach spaces as in 2.8.  By a map between
the invariants $Inv^{0}(A)$ and $Inv^{0}(B)$, we mean a map $$\alpha: (\underline{K}(A),\underline{K}(A)^{+},\Sigma A)\longrightarrow(\underline{K}(B),\underline{K}(B)^{+},\Sigma B)$$ as in 2.7, and for each pair $[p]\in\Sigma A$, $[q]\in\Sigma B$ with $\alpha[p]=[q]$, there is an associate unital positive linear  map (which is automatically continuous as pointed out in 2.8) $$\xi^{p,q}: AffT(pAp)\longrightarrow AffT(qBq)$$ which are compatible in the sense of 2.9 (that is, the diagram (2.A) is commutative for  any pair of projections $p_1\leq p_2$).

\noindent\textbf{2.11.}~~Let $[p]\in\Sigma A$ be represented by $p\in A$. Let $\alpha([p])=[q]$ for $q\in B$. Then $\alpha$ induces a map (still denoted by $\alpha$) $\alpha: K_{0}(pAp)\longrightarrow K_{0}(qBq).$ Note that the  natural map $\rho:=\rho_{pAp}:  K_{0}(pAp)\longrightarrow AffT(pAp)$, defined in 2.8,  satisfies
$\rho(K_{0}(pAp)^{+})\subseteq AffT(pAp)_{+}$ and $\rho([p])=\textbf{1}\in AffT(pAp)$. By 1.20 of [Ji-Jiang], the compatibility in 2.9 (diagram (2.$A$) in 2.9) implies that the following diagram commutes:
$$\CD
  K_{0}(pAp)@>\rho>> AffT(pAp) \\
  @V \alpha VV @V \xi^{p,q} VV  \\
 K_{0}(qBq)@>\rho>> AffT(qBq)~~.
\endCD\;\;\;\;\;\;\;\;\;~~~~~~~~~~~~(2.B)$$
For $p=\textbf{1}_{A}$, this compatibility (the commutativity of diagram (2.$B$)) is included as a part of Elliott invariant for unital simple $C^{*}$-algebras. But
this information are contained in our invariant $Inv^{0}(A)$, as pointed out in [Ji-Jiang].

\noindent\textbf{2.12.}~~Let $A$ be a unital $C^{*}$-algebra,
$B\in \cal{HD}$ and $\{p_{i}\}^{n}_{i=1}\subset B$ be mutually orthogonal projections with $\Sigma p_{i}=\textbf{1}_{B}$. Write $B=\oplus_{j=1}^m B^j$ with $B^j$ being either $PM_{\bullet}(C(X))P$ or $M_l(I_k)$, and for any $i=1,2,\cdots, n$ write $p_i=\oplus_{j=1}^m p_i^j$ with $p_i^j\in B^j$, for $j=1,2,\cdots, m$. Note that for all $\tau\in TB^j$, $\tau(p_i^j)=\frac{rank(p_i^j)}{rank(1_{B^j})}$ (see 2.2 for the definition of $rank$ function), which is independent of $\tau \in TB^j$.

Let $\xi_{i}=(\xi_i^1,\xi_i^2,\cdots, \xi_i^m): AffTA\longrightarrow AffT(p_{i}Bp_{i})=\oplus_{j=1}^m AffT(p_i^j B^j p_i^j)$ be unital positive linear maps, then we can define $\xi=(\xi^1,\xi^2,\cdots, \xi^m): AffTA\longrightarrow AffTB=\oplus_{j=1}^m AffTB^j$ as below
$$\xi^j(f)(\tau)=\sum_{\{i:\tau(p_i^j)\not=0\}} \tau(p_{i}^j)\xi_{i}^j(f)(\frac{\tau|_{p_{i}^jB^jp_{i}^j}}{\tau(p_{i}^j)})~~~\mbox{for} ~~f \in AffTA~~\mbox{and}~~\tau \in TB^j.$$
Note that $\frac{\tau|_{p_{i}^jB^jp_{i}^j}}{\tau(p_{i}^j)}\in T(p_{i}^jB^jp_{i}^j)$. So $\xi_{i}^j(f)$ can evaluate at $\frac{\tau|_{p_{i}^jB^jp_{i}^j}}{\tau(p_{i}^j)}$. Since the value of  $\tau(p_i^j)$ is independent of $\tau\in TB^j$, it is straight forward to verify that $\xi^j\in AffTB^j$.
We denote such $\xi$ by
$\oplus\xi_{i}$.  (For the case that $B$ is general stably finite unital simple $C^*$-algebras with mutually orthogonal projections $\{p_i\}$ with sum $\one_B$, this kind of construction can be carried out by using Lemma 6.4 of [Lin5].)

If $\phi_{i}: A\longrightarrow p_{i}Bp_{i}$ are unital homomorphisms and $\phi=\oplus\phi_{i}: A\longrightarrow B$, then
$$(AffT\phi)^j(f)(\tau)=\sum\limits_{\{i:\tau(p_{i}^j)\neq0\}} \tau(p_{i}^j)AffT\phi_{i}^j(f)(\frac{\tau|_{p_{i}^jB^jp_{i}^j}}{\tau(p_{i}^j)}),$$ where $\phi_i^j: A\to p_i^jB^jp_i^j$ is the j-th component of the map of $\phi_i$. That is,  $AffT\phi=\oplus AffT\phi_{i}$. In particular, if
$\|AffT\phi_{i}(f)-\xi_{i}(f)\|<\varepsilon$ for all $i$, then
$$\|AffT\phi(f)-\xi(f)\|<\varepsilon.$$

\noindent\textbf{2.13.}~~Now, we will introduce the new ingredient of our invariant, which is a simplified version of $U(pAp)/\overline{DU(pAp)}$ for any $[p]\in\Sigma A$, where $DU(pAp)$ is the commutator subgroup of $U(pAp)$. Some notations and prelimary results are quoted from [Thm2], [Thm4] and [NT].

\noindent\textbf{2.14.}~~Let $A$ be a unital $C^{*}$-algebra. Let $U(A)$ denote the group of unitaries of $A$ and,  $U_{0}(A)$, the connected component of $\textbf{1}_{A}$ in $U(A)$. Let $DU(A)$ and $DU_{0}(A)$ denote the commutator subgroups of $U(A)$ and $U_{0}(A)$, respectively. (Recall that the commutator subgroup of a group $G$ is the subgroup generated by all elements of the form $aba^{-1}b^{-1}$, where $a,b\in G$.) One can introduce the following metric
$D_{A}$ on $U(A)/\overline{DU(A)}$ (see [NT,\S3]). For $u,v\in U(A)/\overline{DU(A)}$
$$D_{A}(u,v)=inf\{\|uv^{*}-c\|:~c\in \overline{DU(A)}\},$$
where, on the right hand side of the equation, we use $u,v$ to denote any elements in $U(A)$, which represent the elements $u,v\in U(A)/\overline{DU(A)}$.

\noindent\textbf{Remark 2.15.}~~ Obviously, $D_{A}(u,v)\leq2$. Also,  if $u,v\in U(A)/\overline{DU(A)}$ define two different elements in $K_{1}(A)$, then
$D_{A}(u,v)=2$. (This fact follows from the fact that $\|u-v\|<2$ implies $uv^{*}\in U_{0}(A)$.)

\noindent\textbf{2.16.}~~Let $A$ be a unital $C^{*}$-algebra. Let $AffTA$  and $\rho_A: K_{0}(A)\longrightarrow AffTA$ be as defined as in 2.8,.

For simplicity, we will use $\rho K_0(A)$ to denote the set $\rho_A(K_0(A))$. The metric $d_{A}$ on $AffTA/\overline{\rho K_{0}(A)}$ is defined as follows (see [NT, \S3]).

Let $d^{\prime}$ denote the quotient metric on $AffTA/\overline{\rho K_{0}(A)}$, i.e, for $f,g\in AffTA/\overline{\rho K_{0}(A)}$,
$$d^{\prime}(f,g)=inf\{\|f-g-h\|,h\in\overline{\rho K_{0}(A)}\}.$$
Define $d_{A}$ by
$$d_{A}(f,g)=
 \left\{ \begin{lgathered}
  2,\;\;\; \mbox{if}\;d^{\prime}(f,g)\geq\frac{1}{2} \\
   |e^{2\pi id^{\prime}(f,g)}-1|,\;\;\;\mbox{if}\;d^{\prime}(f,g)<\frac{1}{2}~~~~.
   \end{lgathered} \right.$$
Obviously, $d_{A}(f,g)\leq2\pi d^{\prime}(f,g)$.

\noindent\textbf{2.17.}~~For $A=PM_{k}(C(X))P$, define $SU(A)$ to be the set of unitaries $u\in PM_{k}(C(X))P$ such that for each $x\in X$, $u(x)\in P(x)M_{k}(\mathbb{C})P(x)\cong M_{rank(P)}(\mathbb{C})$ has determinant 1 (note that the determinant of $u(x)$ does not depend on the identification
of $P(x)M_{k}(\mathbb{C})P(x)\cong M_{rank(P)}(\mathbb{C})$). For $A=M_{l}(I_{k})$, by $u\in SU(A)$ we mean that $u\in SU(M_{lk}(C[0,1]))$, where we consider
$A$ to be a subalgebra of $M_{lk}(C[0,1])$. For all basic building blocks  $A\not=M_{l}(I_{k})$, we have $SU(A)=\overline{DU(A)}$. But for $A=M_{l}(I_{k})$, this is not true (see 2.18 and 2.19 below).

In [EGL1], the authors also defined $SU(A)$ for $A$ being a homogeneous algebra and a certain $AH$ inductive limit $C^*$-algebra. This definition can not be generalized to a more general class of $C^*$-algebras. But we will define $\widetilde{SU(A)}$ for any unital $C^*$ algebra $A$. Later,  in our definition of $Inv(A)$, we will only make use of $\widetilde{SU(A)}$ (rather than $SU(A)$).

\noindent\textbf{2.18.}~~Let $A=I_{k}$. Then $K_{1}(A)=\mathbb{Z}/k\mathbb{Z}$, which is generated by $[u]$,  where $u$ is the following unitary
$$u=\left(
      \begin{array}{cccc}
        e^{2\pi i\frac{k-1}{k}t} &  &  &  \\
         & e^{2\pi i(\frac{-t}{k})} &  &  \\
         &  & \ddots &  \\
         &  &  & e^{2\pi i(\frac{-t}{k})}\\
      \end{array}
    \right)\in I_{k}.$$
(Note that $u(0)=\textbf{1}_{k}$, $u(1)=e^{2\pi i(\frac{-1}{k})}\cdot \textbf{1}_{k}$.)

Note that the above $u$ is in $SU(A)$, but not in $U_{0}(A)$,  and therefore not in $DU(A)$.

\noindent\textbf{2.19.}~~By [Thm4] (or [GLN]),  $u\in M_{l}(I_{k})$ is in $\overline{DU(A)}$ if and only if for any irreducible representation $\pi: M_{l}(I_{k})\longrightarrow B(H)$ (dim $H<+\infty$), det$(\pi(u))=1$. For the  unitary $u$ in 2.18, and irreducible representation $\pi$ corresponding
to 1, $\pi(u)=e^{2\pi i(\frac{-1}{k})}$ whose determinant is $e^{2\pi i(\frac{-1}{k})}$ which is not 1. By [Thm2, 6.1] one knows that if $A=I_{k}$,
then $$U_{0}(A)\cap SU(A)=\{e^{2\pi i(\frac{j}{k})},j=0,1,\cdots,k-1\}\cdot\overline{DU(A)}.$$
If $A=M_{l}(I_{k})$, then for any $j\in \mathbb{Z}$,  ~~$e^{2\pi i(\frac{j}{l})}\cdot\textbf{\textbf{1}}_{A}\in\overline{DU(A)}.$ Consequently,
$$U_{0}(A)\cap SU(A)=\{e^{2\pi i(\frac{j}{kl})},j=0,1,\cdots,kl-1\}\cdot\overline{DU(A)}.$$

\noindent\textbf{2.20.}~~Let $\mathbb{T}=\{z\in\mathbb{C},|z|=1\}$. Then for any $A\in\mathcal{HD}$, $\mathbb{T}\cdot\overline{DU(A)}\subset U_{0}(A)$. From 2.17 and 2.19, we have either
$SU(A)=\overline{DU(A)}$ or $U_{0}(A)\cap SU(A)\subset \mathbb{T}\cdot\overline{DU(A)}$.

\noindent\textbf{Lemma 2.21.}~~ Let $A=PM_{k}(C(X))P\in\mathcal{HD}$. For any $u,v\in U(A)$, if $uv^{*}\in \mathbb{T}\cdot\overline{DU(A)}$ (in particular if both
$u,v$ are in $ \mathbb{T}\cdot\overline{DU(A)}$), then $D_{A}(u,v)\leq{2\pi}/{rank(P)}~~.$

Let $A=M_{l}(I_{k})$. For any $u,v$, if  $uv^{*}\in \mathbb{T}\cdot\overline{DU(A)}$, then
 $D_{A}(u,v)\leq{2\pi}/{l}~.$
\begin{proof}
 There is $\omega\in \overline{DU(A)}$ such that $uv^{*}=\lambda\omega$ for some $\lambda\in \mathbb{T}$. Choose $\lambda_{0}=e^{2\pi i\frac{j}{rank(P)}}$, $j\in \N$,  such that $|\lambda-\lambda_{0}|<{2\pi}/{rank(P)}$. And $\lambda_{0}\cdot P\in PM_{k}(C(X))P$ has determinant 1 everywhere and
is in $\overline{DU(A)}$. And so does $\lambda_{0}\omega$. Also we have $|uv^{*}-\lambda_{0}\omega|<{2\pi}/{rank(P)}~.$

The case $A=M_{l}(I_{k})$ is similar.
\end{proof}

\noindent\textbf{2.22.}~~ Let $path(U(A))$ denote the set of piecewise smooth paths $\xi:[0,1]\rightarrow U(A)$. Recall that de la Harp-Skandalis determinant
(see [dH-S]) $\Delta: path(U(A))\rightarrow AffTA$ is defined by
$$\Delta(\xi)(\tau)=\frac{1}{2\pi i}\int_{0}^{1}\tau(\frac{d\xi}{dt}\cdot\xi^{\ast})dt.$$
It is proved in [dH-S](see also [Thm4]) that $\Delta$ induces a map
$\Delta^{\circ}:\pi_{1}(U_{0}(A))\rightarrow AffTA.$
For any two paths $\xi_{1},\xi_{2}$ starting at $\xi_{1}(0)=\xi_{2}(0)=1\in A$ and ending at the same unitary $u=\xi_{1}(1)=\xi_{2}(1)$, we have that
$$\Delta(\xi_{1})-\Delta(\xi_{2})=\Delta(\xi_{1}\cdot\xi_{2}^{\ast})\subset\Delta^{\circ}(\pi_{1}(U_{0}(A))).$$
Consequently $\Delta$ induces a map
\begin{center}
  $\overline{\Delta}: U_{0}(A)\rightarrow AffTA/\Delta^{\circ}(\pi_{1}(U_{0}(A))).$ (See [Thm4, section 3].)
\end{center}
Passing to matrix over $A$, we have a map
$\overline{\Delta}_{n}: U_{0}(M_{n}(A))\rightarrow AffTA/\Delta^{\circ}_{n}(\pi_{1}(U_{0}(M_n(A)))).$

If $1\leq m<n$, then $path(U(M_m(A)))$ (and $U_0(M_m(A))$ ) can be embedded into $path(U(M_n(A)))$ (and $U_0(M_n(A))$ ) by sending $u(t)$ to $diag(u(t),1_{n-m})$. From the above definition, and the formula $$\frac{d}{dt}(diag(u(t),1_{n-m})=diag(\frac{d}{dt}(u(t)), 0_{n-m}),$$ one gets $$\overline{\Delta}_{n}|_{U_0(M_m(A))}=\overline{\Delta}_m.$$
Recall that the Bott isomorphism
$b: K_{0}(A)\rightarrow K_{1}(SA)$
is given by the following:  for any $x\in K_{0}(A)$ represented by a projection $p\in M_{n}(A)$, we have
$$b(x)=[e^{2\pi it}p+(\textbf{1}_{n}-p)]\in K_{1}(SA).$$
If $\xi(t)=e^{2\pi it}p+(\textbf{1}_{n}-p)$, then
$$ (\Delta^{\circ}\xi)(\tau)=\frac{1}{2\pi i}\int_{0}^{1}\tau((2\pi ie^{2\pi it}p)\cdot (e^{-2\pi it}p+(1-p)))dt=\frac{1}{2\pi i}\int_{0}^{1}\tau(2\pi ip)dt=\tau(p).$$
Since Bott map is an isomorphism, it follows  that each loop in $\pi_{1}(U_{0}(A))$ is homotopic to a product of loops of the above form $\xi(t)$.
Consequently $\Delta^{\circ}(\pi_{1}(U_{0}(M_{n}(A))))\subset\rho_{A} K_{0}(A)$. Hence $\overline{\Delta}_{n}$ can be regarded as a map
$$\overline{\Delta}_{n}: U_{0}(M_{n}(A))\rightarrow AffTA/\overline{\rho_{A} K_{0}(A)}~~.$$

\noindent\textbf{Proposition 2.23.}~~  For $A\in\mathcal{HD}$ or $A\in A\mathcal{HD}$, $\overline{DU_{0}(A)}=\overline{DU(A)}$.
\begin{proof}Let the determinant function $\overline{\Delta}_{n}: U_{0}(M_{n}(A))\longrightarrow AffTA/\overline{\Delta^{0}_{n}(\pi_{1}U_{0}(M_{n}(A)))}$ be defined as in $\S3$ of [Thm4] (see 2.22 above). As observed in [NT] (see top of page 33 of [NT]), Lemma 3.1
of [Thm4] implies that $\overline{DU_{0}(A)}=U_{0}(A)\cap\overline{DU(A)}$. For reader's convenience, we give a brief proof of this fact. Namely, the equation $$\left(
                         \begin{array}{ccc}
                           uvu^{-1}v^{-1} & 0 & 0 \\
                           0 & \textbf{1} & 0 \\
                           0 & 0 & \textbf{1} \\
                         \end{array}
                       \right)=\left(
                                 \begin{array}{ccc}
                                   u & 0 & 0 \\
                                   0 & u^{-1} & 0 \\
                                   0 & 0 & \textbf{1} \\
                                 \end{array}
                               \right)\left(
                                        \begin{array}{ccc}
                                          v & 0 & 0 \\
                                          0 & \textbf{1} & 0 \\
                                          0 & 0 & v^{-1} \\
                                        \end{array}
                                      \right)\left(
                                               \begin{array}{ccc}
                                                 u^{-1} & 0 & 0 \\
                                                 0 & u & 0 \\
                                                 0 & 0 & \textbf{1} \\
                                               \end{array}
                                             \right)\left(
                                                      \begin{array}{ccc}
                                                        v^{-1} & 0 & 0 \\
                                                        0 & \textbf{1} & 0 \\
                                                        0 & 0 & v \\
                                                      \end{array}
                                                    \right)$$
implies that $\overline{DU(A)}\subset\overline{DU_{0}(M_{3}(A))}.$ Therefore by Lemma 3.1 of [Thm4], $\overline{DU(A)}\subset ker\overline{\Delta}_{3}.$ If $x\in U_{0}(A)\cap\overline{DU(A)}$, then $\overline{\Delta}_{1}$ is defined at $x$. By calculation in 2.22,  $\overline{\Delta}_{3}|_{U_{0}(A)}=\overline{\Delta}_{1}$. Hence we have $\overline{\Delta}_{1}(x)=0.$ And therefore $x\in\overline{DU_{0}(A)}=ker\overline{\Delta}_{1},$ by Lemma 3.1 of [Thm4].
Note that if $A\in\mathcal{HD}$ or $A\mathcal{HD}$, then  $\overline{DU(A)}\subset U_{0}(A)$.\\
\end{proof}
(It is not known to the authors whether it is always true that $\overline{DU_{0}(A)}=\overline{DU(A)}.$)

\noindent\textbf{2.24.}~~ There is a natural map $\alpha: \pi_{1}(U(A))\longrightarrow K_{0}(A)$, or more generally,
$\alpha_{n}: \pi_{1}(U(M_{n}(A))\longrightarrow K_{0}(A))$ {for any} $n\in \N .$
We need the following notation. For a unital $C^{*}$-algebra $A$, let $\mathcal{P}_{n}K_{0}(A)$ (see [GLX]) be the subgroup of $K_{0}(A)$ generated by the formal
difference of projections $p,q\in M_{n}(A)$ (instead of $M_{\infty}(A)$). Then $$\mathcal{P}_{n}K_{0}(A)\subset Image(\alpha_{n}).$$
In particular,  if $\rho: K_{0}(A)\longrightarrow AffTA$ satisfies $\rho(\mathcal{P}_{n}K_{0}(A))=\rho K_{0}(A)$, then by Theorem 3.2 of [Thm4],
$$U_{0}(M_{n}(A))/\overline{DU_{0}(M_{n}(A))}\cong U_{0}(M_{\infty}(A))/\overline{DU_{0}(M_{\infty}(A))}\cong AffTA/\overline{\rho K_{0}(A)}.$$
Note that for all $A\in\mathcal{HD}$, we have $\rho(\mathcal{P}_{1}K_{0}(A))=\rho K_{0}(A)$ (see below). Consequently, $$U_{0}(A)/\overline{DU_{0}(A)}\cong AffTA/\overline{\rho K_{0}(A)}.$$
If $A$ does not contain building blocks of form $PM_{n}(C(T_{\uppercase\expandafter{\romannumeral2},k}))P$, then such $A$ is the special case of [Thm2], and the above fact is observed in [Thm2] (for circle algebras in [NT] earlier)---in this special case, we ever have $\mathcal{P}_{1}K_{0}(A)=K_{0}(A)$ (as used in [NT] and [Thm2] in the form
of surjectivity of $\alpha: \pi_{1}(U(A))\longrightarrow K_{0}(A)$).
For $A=PM_{n}(C(T_{\uppercase\expandafter{\romannumeral2},k}))P$, we do not have the surjectivity of $\alpha: \pi_{1}(U(A))\longrightarrow K_{0}(A)$ any more. But $K_{0}(A)=\mathbb{Z}\oplus\mathbb{Z}/k\mathbb{Z}$ and image$(\alpha)=\mathcal{P}_{1}K_{0}(A)$ contains at least one element which corresponds to a rank one projection
(any bundle over $T_{\uppercase\expandafter{\romannumeral2},k}$ has a subbundle of rank 1)---that is, $$\rho(\mathcal{P}_{1}K_{0}(A))=\rho K_{0}(A)(\subseteq AffTA)$$ consisting all constant functions
from $T_{\uppercase\expandafter{\romannumeral2},k}$ to $\frac{1}{rank(P)}\mathbb{Z}$.

As in [NT, Lemma 3.1] and [Thm 2, Lemma 6.4], the map $\overline{\Delta}:~ U_{0}(A)\rightarrow AffTA/\overline{\rho_{A}(K_{0}(A))}$ (in 2.22) has Ker$\overline{\Delta}=\overline{DU(A)}$ and the following lemma holds.

\noindent\textbf{Lemma 2.25.}~~ If a unital $C^{*}$-algebra $A$ satisfies $\rho(\mathcal{P}_{1}K_{0}(A))=\rho K_{0}(A)$ and $\overline{DU_{0}(A)}=\overline{DU(A)}$ (see 2.24 and 2.23), in particular, if $A\in\mathcal{HD}$ or $A\in A\mathcal{HD}$, then the following hold:\\
(1)~~ There is a split exact sequence $$0\rightarrow AffTA/\overline{\rho K_{0}(A)}\xrightarrow{\lambda_{A}} U(A)/\overline{DU(A)}\rightarrow K_{1}(A)\rightarrow0.$$
(2)~~ $\lambda_{A}$ is an isometry with respect to the metrics $d_{A}$ and $D_{A}$.

\noindent\textbf{2.26.}~~ Recall from $\S3$ of [Thm4], the de la Harpe---Skandalis determinant (see [dH-S]) can be used to define $$\overline{\Delta}: U_{0}(A)/\overline{DU(A)}\longrightarrow AffTA/\overline{\rho K_{0}(A)}.$$ With the condition of  Lemma 2.25 above, this map is an isometry with respect to the metrics $d_{A}$ and $D_{A}$.  In fact, the inverse of this map is $\lambda_{A}$ in  Lemma 2.25.

It follows from the definition of $\overline{\Delta}$ (see $\S3$ of [Thm4]) that
$$\overline{\Delta}(e^{2\pi itp})=t\cdot\rho([p])~~~~~~~\;\;(mod~(\overline{\rho K_{0}(A))}).\;\;\;\;~~~~~~~~~~~~~~~~~~~~~~~~~~~~~~~~~~(2.c)$$
where $[p]\in K_{0}(A)$ is the element represented by projection
$p\in A$.

It is convenient to introduce the extended commutator group $DU^{+}(A)$, which is generated by \\$DU(A)\subset U(A)$ and the set
$\{e^{2\pi itp}=e^{2\pi it}p+(\textbf{1}-p)\in U(A)~|~~t\in\mathbb{R},p\in A ~\mbox{ is a projection}\}.$
Let $\widetilde{DU(A)}$ denote the closure of $DU^{+}(A)$. That is,  $\widetilde{DU(A)}=\overline{DU^{+}(A)}$.

Let us use $\widetilde{\rho K_{0}(A)}$
to denote the real vector space spanned by $\overline{\rho K_{0}(A)}$. That is,  $$\widetilde{\rho K_{0}(A)}:=\overline{\{\Sigma\lambda_{i}\phi_{i}~|~~\lambda_{i}\in\mathbb{R},\phi_{i}\in\rho K_{0}(A)}\}.$$

Suppose that $\overline{\rho K_{0}(A)}=\overline{\rho(\mathcal{P}_{1}K_{0}(A))}$. It follows from (2.c), the image of $\widetilde{DU(A)}/\overline{DU(A)}$ under the map $\overline{\Delta}$ is
exactly
$\widetilde{\rho K_{0}(A)}/\overline{\rho K_{0}(A)}$. Therefore $\lambda_A$ takes $\widetilde{\rho K_{0}(A)}/\overline{\rho K_{0}(A)}$ to $\widetilde{DU(A)}/\overline{DU(A)}$.
Hence
$\overline{\Delta}:U_0(A)/\overline{DU(A)}\longrightarrow AffTA/\overline{\rho K_{0}(A)}$ also induces a quotient map (still denoted by $\overline{\Delta}$)
$$\overline{\Delta}: U_{0}(A)/\widetilde{DU(A)}\longrightarrow AffTA/\widetilde{\rho K_{0}(A)}$$ which is an isometry using the quotient metrics of $d_{A}$ and $D_{A}$. The inverse of this quotient map $\overline{\Delta}$ gives rise to the isometry
$$\widetilde{\lambda}_{A}: AffTA/\widetilde{\rho K_{0}(A)}\longrightarrow U_{0}(A)/\widetilde{DU(A)}\hookrightarrow U(A)/\widetilde{DU(A)}$$
which is an isometry with respect to the quotient metrics $\widetilde{d}_{A}$ and $\overline{D_{A}}$ as described below.

For any $u,v\in U(A)/\widetilde{DU(A)},$ $$\overline{D_{A}}(u,v)=inf\{\|uv^{*}-c\|~|~~c\in\widetilde{DU(A)}\}.$$ Let $\widetilde{d'}$ denote the
quotient metric on $AffTA/\widetilde{\rho K_{0}(A)}$ of $AffTA$, that is,
$$\widetilde{d'}(f,g)=inf\{\|f-g-h\|~|~~h\in\widetilde{\rho K_{0}(A)}\}~~~~~\forall f,g\in AffTA/\widetilde{\rho K_{0}(A)}.$$ Define $\widetilde{d}_{A}$ by
$$\widetilde{d}_{A}(f,g)=
 \left\{ \begin{lgathered}
  2,\;\;\; \;\;\;\;\;\;\;\;\;\;\;\;\;\;\;\;\;\;if\;\widetilde{d'}(f,g)\geq\frac{1}{2} \\
  |e^{2\pi i\widetilde{d'}(f,g)}-1|,\;\;\;if\;\widetilde{d'}(f,g)<\frac{1}{2}~~~.
   \end{lgathered} \right.$$

The following result is a consequence of Lemma 2.25.\\
\noindent\textbf{Lemma 2.27.}~~ If a unital $C^{*}$-algebra $A$ satisfies $\rho(\mathcal{P}_{1}K_{0}(A))=\rho K_{0}(A)$ and $\overline{DU_{0}(A)}=\overline{DU(A)}$ (see 2.24 and 2.23), in particular, if $A\in\mathcal{HD}$ or $A\in A\mathcal{HD}$, then  we have\\
(1)~~ There is a split exact sequence $$0\rightarrow AffTA/\widetilde{\rho K_{0}(A)}\xrightarrow {\widetilde{\lambda}_{A}}U(A)/\widetilde{DU(A)}\xrightarrow{\pi_{A}} K_{1}(A)\rightarrow0.$$
(2)~~ $\widetilde{\lambda}_{A}$ is an isometry with respect to $\widetilde{d}_{A}$ and $\overline{D_{A}}$.

\begin{proof} As we mentioned in 2.26, the map $\lambda_A$ in Lemma 2.25 takes $\widetilde{\rho K_{0}(A)}/\overline{\rho K_{0}(A)}$ to $\widetilde{DU(A)}/\overline{DU(A)}$. From the exact sequence in Lemma 2.25, passing  to quotient, one gets the exact sequence in (1).

Note that $\widetilde{d}_{A}$ on $AffTA/\widetilde{\rho K_{0}(A)}$ is the quotient metric induced by $d_A$ on $AffTA/ \overline{\rho K_{0}(A)}$ and $\overline{D_{A}}$ on $U(A)/\widetilde{DU(A)}$ is the quotient metric induced by $D_A$ on $U(A)/\overline{DU(A)}$. Hence $\widetilde{\lambda}_{A}$ is an isometry, since so is $\lambda_A$.

\end{proof}

\noindent\textbf{2.28.}~~ Instead of $\widetilde{DU(A)}$, we will need the group
$$\widetilde{SU(A)}:=\{\overline{x\in U(A)~|~~x^{n}\in\widetilde{DU(A)}~~~~~\mbox{for some} ~n\in\mathbb{Z}_{+}\backslash\{0\}}\}.$$

For $A\in\mathcal{HD}$, say $A=PM_{l}(C(X))P$  ($X=[0,1], S^1$ or $T_{II,k}$) or  $A=M_{l}(I_{k})$, ~~$\widetilde{SU(A)}$ is the set of all unitaries $u\in P(M_{l}C(X))P$ or $u\in M_{l}(I_{k})$ such that the determinant function $$X\ni x\longmapsto det(u(x))~~~~~~\mbox{or}~~~~~(0,1)\ni t\longmapsto det(u(t))$$
is a constant function.
Comparing with the set $SU(A)$ in [EGL1] or 2.17 above (which only defines for $\mathcal{HD}$ blocks), where the function will be constant 1, here we allow the function to be arbitrary constant in $\mathbb{T}$. Hence  for a basic building block  $A=PM_n(C(X))P\in \cal{HD}$ or $A=M_l(I_k)$,   $$\widetilde{SU(A)}=\mathbb{T}\cdot SU(A).$$

The notations $\widetilde{\rho K_0(A)}$, $\widetilde{DU(A)}$ and $\widetilde{SU(A)}$ reflect that they are constructed from $\rho K_0(A)$, $DU(A)$ and $SU(A)$, respectively.
{\bf To make the notation simpler, from now on, we will use $\widetilde{\rho K_0}(A)$ to denote $\widetilde{\rho K_0(A)}=\widetilde{\rho_A(K_0(A))}$, $\widetilde{DU}(A)$ to denote $\widetilde{DU(A)}$, and $\widetilde{SU}(A)$ to denote $\widetilde{SU(A)}$.}

\noindent\textbf{Lemma 2.29.}~~ Let $\alpha,\beta: K_{1}(A)\longrightarrow U(A)/\widetilde{DU}(A)$ be two splittings of $\pi_{A}$ in Lemma 2.27. Then $$\alpha|_{tor\;K_{1}(A)}=\beta|_{tor\;K_{1}(A)}$$ and $\alpha(tor\;K_{1}(A))\subset\widetilde{SU}(A)/\widetilde{DU}(A)$. Furthermore,  $\alpha$ identifies $tor(K_{1}(A))$ with $\widetilde{SU}(A)/\widetilde{DU}(A)$.
\begin{proof}
For any $z\in tor\;K_{1}(A)$, with $kz=0$ for some integer $k>0$, we have $$\pi_{A}\alpha(z)=z=\pi_{A}\beta(z).$$
By the exactness of the sequence, there is an element $f\in AffTA/\widetilde{\rho K_{0}}(A)$ such that $$\alpha(z)-\beta(z)=\widetilde{\lambda}_{A}(f).$$ Since $k\alpha(z)-k\beta(z)=\alpha(kz)-\beta(kz)=0,$ we have $\widetilde{\lambda}_{A}(kf)=0$. By the injectivity of $\widetilde{\lambda}_{A}$, $kf=0$. Note that
$AffTA/\widetilde{\rho K_{0}}(A)$ is an $\mathbb{R}$-vector space, $f=0$. Furthermore, $k\alpha(z)=0$ in $U(A)/\widetilde{DU}(A)$ implies that
$$\alpha(z)\in\widetilde{SU}(A)/\widetilde{DU}(A).$$ Hence we get $\alpha(tor\;K_{1}(A))\subset\widetilde{SU}(A)$. If $u\in\widetilde{SU}(A)/\widetilde{DU}(A)$ then $\alpha(\pi_{A}(u))=u$.\\
\end{proof}

\noindent\textbf{2.30.}~~ Let $U_{tor}(A)$ denote the set of unitaries $u\in A$ such that $[u]\in tor\;K_{1}(A)$.
For any $C^*$ algebra $A$ we have
$\widetilde{SU}(A)\subset U_{tor}(A).$ If we further assume $\overline{DU_{0}(A)}=\overline{DU(A)}$, then$$ ~~\widetilde{DU}(A)=U_{0}(A)\cap\widetilde{SU}(A)~~\mbox{ and}~~U_{tor}(A)=U_{0}(A)\cdot\widetilde{SU}(A).$$ Evidently, we have $U_{0}(A)/\widetilde{DU}(A)\cong U_{tor}(A)/\widetilde{SU}(A).$
The metric $\overline{D_{A}}$ on $U(A)/\widetilde{DU}(A)$ induces a metric $\widetilde{D}_{A}$ on $U(A)/\widetilde{SU}(A)$. And the above identification $U_{0}(A)/\widetilde{DU}(A)$ with $U_{tor}(A)/\widetilde{SU}(A)$ is an isometry with respect to $\overline{D}_A$ and $\widetilde{D}_{A}$.
Hence $\widetilde{\lambda}_{A}$ in 2.26 can be regarded as a map (still denoted  by $\widetilde{\lambda}_{A}$):
$$\widetilde{\lambda}_{A}: AffTA/\widetilde{\rho K_{0}}(A)\longrightarrow U_{tor}(A)/\widetilde{SU(A)}\hookrightarrow U(A)/\widetilde{SU(A)}.$$

Similar to Lemma 2.27, we have

\noindent\textbf{Lemma 2.31.}~~ If a unital $C^{*}$-algebra $A$ satisfies $\rho(\mathcal{P}_{1}K_{0}(A))=\rho K_{0}(A)$ and $\overline{DU_{0}(A)}=\overline{DU(A)}$ (see 2.24 and 2.23), in particular, if $A\in\mathcal{HD}$ or $A\in A\mathcal{HD}$, then the following hold:\\
(1)~~ There is a split exact sequence $$0\rightarrow AffTA/\widetilde{\rho K_{0}}(A)\xrightarrow{\widetilde{\lambda}_{A}}U(A)/\widetilde{SU}(A)\xrightarrow{\pi_{A}} K_{1}(A)/tor\;K_{1}(A)\rightarrow0.$$
(2)~~ $\widetilde{\lambda}_{A}$ is an isometry with respect to the metrics $\widetilde{d}_{A}$ and $\widetilde{D}_{A}$.

\noindent\textbf{2.32.}~~ For each pair of projections $p_{1},p_{2}\in A$ with $p_{1}=up_{2}u^{*},$ $$U(p_{1}Ap_{1})/\widetilde{SU}(p_{1}Ap_{1})\cong U(p_{2}Ap_{2})/ \widetilde{SU}(p_{2}Ap_{2}).$$ Also,  since in any unital $C^{*}$-algebra $A$ and unitaries $u,v\in U(A)$, $v$ and $uvu^{*}$ represent
 a same element in $U(A)/\widetilde{SU}(A)$, and the above identification does not depend on the choice
of $u$ to implement $p_{1}=up_{2}u^{*}$. That is for any $[p]\in\Sigma A$, the group $U(pAp)/\widetilde{SU}(pAp)$ is well defined, which does not depend on choice of $p\in[p]$. We will include this group (with metric) as part of our invariant. If $[p]\leq[q]$, then we can choose $p,q$ such that $p\leq q$. In this case,  there is a natural inclusion map $\imath: pAp\longrightarrow qAq$  which induces $$\imath_{*}: U(pAp)/\widetilde{SU}(pAp)\longrightarrow U(qAq)/\widetilde{SU}(qAq),$$ where $\imath_{*}$ is defined by
$$\imath_{*}(u)=u\oplus(q-p)\in U(qAq),~~~~ \forall u\in U(pAp).$$
~~~A unital homomorphism $\phi: A\to B$ induces a contractive group homomorphism
$$\phi^{\natural}: U(A)/\widetilde{SU}(A) \longrightarrow U(B)/\widetilde{SU}(B).$$
If $\phi$ is not  unital, then the map $\phi^{\natural}: U(A)/\widetilde{SU}(A) \longrightarrow U(\phi(\textbf{1}_A)B\phi(\textbf{1}_A))/\widetilde{SU}(\phi(\textbf{1}_A)B\phi(\textbf{1}_A))$ is induced by the corresponding unital homomorphism.  In this case,  $\phi$ also induces the map $\imath_*\circ \phi^{\natural}: U(A)/\widetilde{SU}(A) \longrightarrow U(B)/\widetilde{SU}(B)$, which is denoted by $\phi_*$ to avoid confusion. If $\phi$ is unital, then $\phi^{\natural}=\phi_*$. If $\phi$ is not unital, then  $\phi^{\natural}$ and $\phi_*$ have different codomains.
 That is, $\phi^{\natural}$ has codomain $U(\phi(\textbf{1}_A)B\phi(\textbf{1}_A))/\widetilde{SU}
(\phi(\textbf{1}_A)B\phi(\textbf{1}_A))$, but $\phi_*$ has codomain $U(B)/\widetilde{SU}(B)$. (See some further explanation with an example in the last paragraph of  3.7 below.)

Since $U(A)/\widetilde{SU}(A) $ is an Abelian group, we will call the unit $[{\bf 1}]\in U(A)/\widetilde{SU}(A) $  the zero element. If $\phi: A \to B$ satisfies $\phi(U(A))\subset \widetilde{SU}(\phi(\textbf{1}_A)B\phi(\textbf{1}_A))$, then $\phi^{\natural}=0$. In particular, if the image of $\phi$ is of finite dimensional, then $\phi^{\natural}=0$.

\noindent\textbf{2.33.}~~ In this paper and [GJL], we will denote $$(\underline{K}(A),\underline{K}(A)^{+},\Sigma A,\{AffT(pAp)\}_{[p]\in\Sigma A},\{U(pAp)/\widetilde{SU}(pAp)\}_{[p]\in\Sigma A})$$ by $Inv(A)$. By a map from $Inv(A)$ to $Inv(B)$, we mean
$$\alpha: (\underline{K}(A),\underline{K}(A)^{+},\Sigma A)\longrightarrow(\underline{K}(B),\underline{K}(B)^{+},\Sigma B)$$ as in 2.7, and for each pair $([p], [\overline{p}]) \in\Sigma A\times \Sigma B$ with $\alpha([p])=[\overline{p}]$, there exist an associate unital positive (continuous) linear map
$$\xi^{p,\overline{p}}: AffT(pAp)\longrightarrow AffT(\overline{p}B\overline{p})$$ and an associate contractive group homomorphism
$$\chi^{p,\overline{p}}: U(pAp)/\widetilde{SU}(pAp)\longrightarrow U(\overline{p}B\overline{p})/\widetilde{SU}(\overline{p}B\overline{p})$$
satisfying  the following compatibility conditions. (Note that $\chi^{p,\overline{p}}$ is continuous, as it is a contractive group homomorphism  from a metric group to another metric group.) \\
(a)~~ If $p<q$, then the diagrams
$$\CD
  AffT(pAp) @>\xi^{p,\overline{p}}>> AffT(\overline{p}B\overline{p}) \\
  @V\imath_T  VV @V \imath_T VV \\
  AffT(qAq) @>\xi^{q,\overline{q}}>> AffT(\overline{q}B\overline{q})
\endCD\;\;\;\;\;\;\;\;\;\;\;\;\;\;\;\;\;\;\;\;\;\;\;\;\;~~~~~~~~~(\uppercase\expandafter{\romannumeral1})$$ and
$$\CD
  U(pAp)/\widetilde{SU}(pAp) @>\chi^{p,\overline{p}}>> U(\overline{p}B\overline{p})/\widetilde{SU}(\overline{p}B\overline{p}) \\
  @V\imath_*  VV @V \imath_* VV   \\
  U(qAq)/\widetilde{SU}(qAq) @>\chi^{q,\overline{q}}>> U(\overline{q}B\overline{q})/\widetilde{SU}(\overline{q}B\overline{q})
\endCD\;\;\;\;\;\;\;~~~~~~~~~~~(\uppercase\expandafter{\romannumeral2})$$
commutes, where the vertical maps are induced by inclusions.\\
(b)~~ The following diagram commutes
$$\CD
  K_{0}(pAp) @>\rho>> AffT(pAp) \\
  @V \alpha VV @V \xi^{p,\overline{p}} VV  \\
  K_{0}(\overline{p}B\overline{p}) @>\rho>> AffT(\overline{p}B\overline{p})
\endCD\;\;\;\;\;\;\;\;\;\;\;\;\;\;\;\;\;\;\;\;\;\;\;\;\;\;\;\;\;\;\;\;\;\;~~~~~~~(\uppercase\expandafter{\romannumeral3})$$
and therefore $\xi^{p,\overline{p}}$ induces a map (still denoted by $\xi^{p,\overline{p}}$):
$$\xi^{p,\overline{p}}: AffT(pAp)/\widetilde{\rho K_{0}}(pAp)\longrightarrow AffT(\overline{p}B\overline{p})/\widetilde{\rho K_{0}}(\overline{p}B\overline{p}).$$
(The commutativity  of $(\uppercase\expandafter{\romannumeral3})$ follows from the commutativity of $(\uppercase\expandafter{\romannumeral1})$, by 1.20 of [Ji-Jiang]. So this is not an extra requirement.)\\
(c)~~ The following diagrams
$$\CD
  AffT(pAp)/\widetilde{\rho K_{0}}(pAp) @>>> U(pAp)/\widetilde{SU}(pAp) \\
  @V \xi^{p,\overline{p}} VV @V \chi^{p,\overline{p}} VV  \\
  AffT(\overline{p}B\overline{p})/\widetilde{\rho K_{0}}(\overline{p}B\overline{p}) @>>> U(\overline{p}B\overline{p})/\widetilde{SU}(\overline{p}B\overline{p})
\endCD\;\;\;\;\;~~~~~~~~~~~(\uppercase\expandafter{\romannumeral4})$$ and
$$\CD
  U(pAp)/\widetilde{SU}(pAp) @>>> K_{1}(pAp)/tor\;K_{1}(pAp) \\
  @V\chi^{p,\overline{p}} VV @V \alpha_{1} VV  \\
  U(\overline{p}B\overline{p})/\widetilde{SU}(\overline{p}B\overline{p}) @>>> K_{1}(\overline{p}B\overline{p})/tor\;K_{1}(\overline{p}B\overline{p})
\endCD\;\;\;\;\;\;\;\;\;~~~~~~~~~~~(\uppercase\expandafter{\romannumeral5})$$
commute, where $\alpha_{1}$ is induced by $\alpha$.

We will denote the map from $Inv(A)$ to $Inv(B)$ by
$$(\alpha,\xi,\chi): (\underline{K}(A),\{AffT(pAp)\}_{[p]\in\Sigma A},\{U(pAp)/\widetilde{SU}(pAp)\}_{[p]\in\Sigma A})\longrightarrow~~~~~~~~~~~~~~~~~~~~$$
$$~~~~~~~~~~~~~(\underline{K}(B),\{AffT(\overline{p}B\overline{p})\}_{[\overline{p}]\in\Sigma B},\{U(\overline{p}B\overline{p})/\widetilde{SU}(\overline{p}B\overline{p})\}_{[\overline{p}]\in\Sigma B}).$$
Completely similar to [NT, Lemma 3.2] and [Thm2, Lemma 6.5], we have the following propositions.

\noindent\textbf{Proposition 2.34.}~~ Let  unital $C^{*}$-algebra $A$  ($B$, resp.) satisfy $\rho(\mathcal{P}_{1}K_{0}(A))=\rho K_{0}(A)$ ($\rho(\mathcal{P}_{1}K_{0}(B))=\rho K_{0}(B)$, resp.) and $\overline{DU_{0}(A)}=\overline{DU(A)}$ ($\overline{DU_{0}(B)}=\overline{DU(B)}$, resp.). In particular, let $A,B\in\mathcal{HD}$ or $A\mathcal{HD}$ be unital $C^{*}$-algebras. Assume that $$\psi_{1}: K_{1}(A)\longrightarrow K_{1}(B)~~~\mbox{and }~~~~\psi_{0}: AffTA/\overline{\rho K_{0}(A)}\longrightarrow AffTB/\overline{\rho K_{0}(B)}$$ are group homomorphisms such that $\psi_{0}$ is a
contraction with respect to $d_{A}$ and $d_{B}$. Then there is a group homomorphism $$\psi: U(A)/\overline{DU(A)}\longrightarrow U(B)/\overline{DU(B)}$$ which is a contraction with respect to $D_{A}$ and $D_{B}$ such that the diagram
$$
\xymatrix{
      0\ar[r] & AffTA/\overline{\rho K_{0}(A)}\ar[d]^{\psi_{0}}\ar[r]^-{\lambda_{A}}& U(A)/\overline{DU(A)}\ar[d]^{\psi}\ar[r]^-{\pi_{A}} & K_{1}(A)\ar[d]^{\psi_{1}}\ar[r]& 0 \\
      0\ar[r] & AffTA/\overline{\rho K_{0}(B)}\ar[r]^-{\lambda_{B}}& U(B)/\overline{DU(B)}\ar[r]^-{\pi_{B}} & K_{1}(B)\ar[r] & 0 \\
  }
$$
commutes. If $\psi_{0}$ is an isometric isomorphism and $\psi_{1}$ is an isomorphism, then $\psi$ is an isometric isomorphism.

\noindent\textbf{Proposition 2.35.}~~ Let  unital $C^{*}$-algebra $A$  ($B$, resp.) satisfy $\rho(\mathcal{P}_{1}K_{0}(A))=\rho K_{0}(A)$ ($\rho(\mathcal{P}_{1}K_{0}(B))=\rho K_{0}(B)$, resp.) and $\overline{DU_{0}(A)}=\overline{DU(A)}$ ($\overline{DU_{0}(B)}=\overline{DU(B)}$, resp.). In particular, let $A,B\in\mathcal{HD}$ or $A\mathcal{HD}$ be unital $C^{*}$-algebras. Assume that $$\psi_{1}: K_{1}(A)\longrightarrow K_{1}(B)~~~\mbox{ and }~~~\psi_{0}: AffTA/\widetilde{\rho K_{0}}(A)\longrightarrow AffTB/\widetilde{\rho K_{0}}(B)$$ are group homomorphisms such that $\psi_{0}$ is a
contraction with respect to $\widetilde{d}_{A}$ and $\widetilde{d}_{B}$. Then there is a group homomorphism $$\psi: U(A)/\widetilde{SU}(A)\longrightarrow U(B)/\widetilde{SU}(B)$$ which is a contraction with respect to $\widetilde{D}_{A}$ and $\widetilde{D}_{B}$ such that the diagram
$$
\xymatrix{
      0\ar[r] & AffTA/\widetilde{\rho K_{0}}(A)\ar[d]^{\psi_{0}}\ar[r]^{\widetilde{\lambda}_{A}} & U(A)/\widetilde{SU}(A)\ar[d]^{\psi}\ar[r]^{\widetilde{\pi}_{A}}  & K_{1}(A)/tor\;K_{1}(A)\ar[d]^{\psi_{1}}\ar[r]& 0 \\
      0\ar[r] & AffTA/\widetilde{\rho K_{0}}(B)\ar[r]^{\widetilde{\lambda}_{B}}& U(B)/\widetilde{SU}(B)\ar[r]^{\widetilde{\pi}_{B}} & K_{1}(B)/tor\;K_{1}(B)\ar[r] & 0 \\
  }
$$
commutes. If $\psi_{0}$ is an isometric isomorphism and $\psi_{1}$ is an isomorphism, then $\psi$ is an isometric isomorphism.

\noindent\textbf{Remark 2.36.}~~As in Proposition 2.35 (or Proposition 2.34), for each fixed pair $p\in A$, $\overline{p}\in B$ with $$\alpha([p])=[\overline{p}],$$ if we have an isometric isomorphism between $AffT(pAp)/\widetilde{\rho K_{0}}(pAp)$ and $AffT(\overline{p}B\overline{p})/ \widetilde{\rho K_{0}}(\overline{p}B\overline{p})$ (or between $AffT(pAp)/\overline{\rho K_{0}(pAp)}$ and $AffT(\overline{p}B\overline{p})/\overline{\rho K_{0}(\overline{p}B\overline{p})}$) and isomorphism between $K_{1}(pAp)$ and $K_{1}(pBp)$, then we have an
isometric isomorphism between $U(pAp)/\widetilde{SU}(pAp)$ and $U(\overline{p}B\overline{p})/\widetilde{SU}(\overline{p}B\overline{p})$
(or $U(pAp)/  \overline{DU(pAp)}$ and $U(\overline{p}B\overline{p})/\overline{DU(\overline{p}B\overline{p})}$)  making both diagrams~ $(\uppercase\expandafter{\romannumeral4})$ and
$(\uppercase\expandafter{\romannumeral5})$ commute. This is the reason $U(A)/\overline{DU(A)}$ is not included in the Elliott invariant in the classification of simple $C^{*}$-algebras. For our setting, even though for each pair of projections $(p,\bar{p})$ with $\alpha([p])=[\bar{p}]$, we can find an isometric isomorphism between $U(pAp)/\widetilde{SU}(pAp)$ and
$U(\overline{p}B\overline{p})/\widetilde{SU}(\overline{p}B\overline{p})$, provided that the other parts of invariants $Inv^{0}(A)$ and $Inv^{0}(B)$ are
isomorphic, we still can not make such system of isometric isomorphisms compatible---that is, can not  make the diagram $\uppercase\expandafter{\romannumeral2}$ commutes for $p<q$. We will present
two non isomorphic $C^{*}$-algebras $A$ and $B$ in our class such that $Inv^{0}(A)\cong Inv^{0}(B)$, in next section, where  $Inv^{0}(B)$ is defined in 2.10. Hence it is
essential to include $\{U(pAp)/\widetilde{SU}(pAp)\}_{p\in\Sigma}$ with the compatibility as part of $Inv(A)$.

\noindent\textbf{2.37.}~~Replacing $U(pAp)/\widetilde{SU}(pAp)$, one can also use $U(pAp)/\overline{DU(pAp)}$ as the part of the invariant. That is, one can define $Inv'(A)$ as $$(\underline{K}(A),\underline{K}(A)^{+},\Sigma A,\{AffT(pAp)\}_{[p]\in\Sigma A},\{U(pAp)/\overline{DU(pAp)}\}_{[p]\in\Sigma A}),$$
with corresponding compatibility condition---one needs to change diagrams $(IV)$ and $(V)$ to the corresponding ones. It is not difficult to see that
$Inv'(A)\cong Inv'(B)$ implies $Inv(A)\cong Inv(B)$. We choose the formulation of $Inv(A)$, since it is much more convenient for the proof of the main theorem in [GJL] and it is formally a weaker requirement than the one to require the isomorphism between $Inv'(A)$ and $Inv'(B)$, and the theorem is formally stronger.
(Let us point out that, in the construction of the example (and its proof)  in  section 3 of this article, $Inv'(A)$ is as convenient as $Inv(A)$, and therefore if only for the sake of example in section 3 of this paper,  it is not necessary  to introduce $\widetilde{SU}(A)$.)

Furthermore, it is straight forward to check the following proposition:

\noindent\textbf{Proposition 2.38.} Let  unital $C^{*}$-algebra $A$  ($B$, resp.) satisfy $\rho(\mathcal{P}_{1}K_{0}(A))=\rho K_{0}(A)$ ($\rho(\mathcal{P}_{1}K_{0}(B))=\rho K_{0}(B)$, resp.) and $\overline{DU_{0}(A)}=\overline{DU(A)}$ ($\overline{DU_{0}(B)}=\overline{DU(B)}$, resp.). In particular, let $A,B\in\mathcal{HD}$ or $A\mathcal{HD}$ be unital $C^{*}$-algebras. Suppose that $K_{1}(A)=tor(K_{1}(A))$ and $K_{1}(B)=tor(K_{1}(B))$. It follows that
$Inv^0(A)\cong Inv^0(B)$ implies that $Inv(A)\cong Inv(B)$.
\begin{proof}
It follows from the fact that any isomorphism
\begin{center}
$\xi^{p,\overline{p}}:~~AffT(pAp)/\widetilde{\rho K_0}(pAp)\longrightarrow
 AffT(\overline{p}B\overline{p})/\widetilde{\rho K_0}(\overline{p}B\overline{p})$
\end{center}
induces a unique isomorphism
\begin{center}
$\chi^{p,\overline{p}}: U(pAp)/\widetilde{SU}(pAp)\longrightarrow U(\overline{p}B\overline{p})/\widetilde{SU}(\overline{p}B\overline{p})$
\end{center}
(Note that by the split exact sequence in Lemma 2.31, we have $AffT(pAp)/\widetilde{\rho K_0}(pAp)\cong U(pAp)/\widetilde{SU}(pAp)$).\\
\end{proof}

\vspace{0.1in}

The following calculations and notations will be used in [GJL].

\noindent\textbf{2.39.}~~In general, for $A=\oplus A^{i}$, $\widetilde{SU}(A)=\oplus_{i}\widetilde{SU}(A^{i}).$ For $A=PM_{l}(C(X))P\in\mathcal{HD}$, $\widetilde{SU}(A)=\widetilde{DU}(A).$ For $A=M_{l}(I_{k})$, $\widetilde{SU}(A)=\widetilde{DU}(A)\oplus K_{1}(A).$
For both cases, $U(A)/\widetilde{SU}(A)$ can be identified with \\$C_{1}(X,S^{1}):=C(X,S^{1})/\{constant\; functions\}$ or with
$C_{1}([0,1],S^{1})=C([0,1],S^{1})/\{constant\; functions\},$ for $A=M_l(I_{k}).$

Furthermore, $C_{1}(X,S^{1})$ can be identified as the set of continuous functions from $X$ to
$S^{1}$ such that $f(x_{0})=1$ for certain fixed base point $x_{0}\in X$. For $X=[0,1]$, we choose 0 to be the base point. For $X=S^{1}$, we choose $1\in S^{1}$ to be the base point.

\noindent\textbf{2.40.}~~Let $A=\oplus^{n}_{i=1}A^{i}\in\mathcal{HD}$, $B=\oplus^{m}_{j=1}B^{j}\in\mathcal{HD}$. In this subsection we will discuss some consequences of the compatibility of the maps between $AffT$ spaces. Let $$p=\oplus p^{i}<q=\oplus q^{i}\in A~~~\mbox{ and }~~~\overline{p}=\oplus_{j=1}^{m} \overline{p}^{j}<\overline{q}=\oplus_{j=1}^{m} \overline{q}^{j}\in B$$ be projections satisfying $\alpha([p])=[\overline{p}]$ and $\alpha([q])=[\overline{q}]$. Suppose that two unital positive linear maps $\xi_{1}: AffTpAp\longrightarrow AffT\overline{p}B\overline{p}$ and
$\xi_{2}: AffTqAq\longrightarrow AffT\overline{q}B\overline{q}$ are compatible with $\alpha$ (see diagram (2.$B$) in 2.11) and compatible with each
other (see diagram (2.$A$) in 2.9). Since the (not necessarily unital) maps $AffTpAp\longrightarrow AffTqAq$ and $AffT\overline{p}B\overline{p}\longrightarrow AffT\overline{q}B\overline{q}$ induced by inclusions, are injective, we know that the map $\xi_{1}$
is completely decided by $\xi_{2}$. Let $$\xi_{2}^{i,j}: AffTq^{i}Aq^{i}\longrightarrow AffT\overline{q^{j}}B^{j}\overline{q^{j}}~~(\mbox{ or}~~ \xi_{1}^{i,j}: AffTp^{i}Ap^{i}\longrightarrow AffT\overline{p^{j}}B^{j}\overline{p^{j}})$$ be the corresponding component of the map $\xi_{2}$ (or $\xi_{1}$). If $p^{i}\neq0$
and $\overline{p}^{j}\neq0$, then $\xi_{1}^{i,j}$ is given by the following formula, for any $f\in AffTp^{i}A^{i}p^{i}=C_{\mathbb{R}}(SpA^{i})(\cong AffTq^{i}Aq^{i})$,
 $$\xi_{1}^{i,j}(f)=\frac{rank\;\overline{q_{j}}}{rank\;\overline{p_{j}}}\cdot\frac{rank\;\alpha^{i,j}(p^{i})}{rank\;\alpha^{i,j}(q^{i})}\cdot\xi_{2}^{i,j}(f).$$
In particular, if $q=\textbf{1}_{A}$ with $\overline{q}=\alpha_{0}[\textbf{1}_{A}]$,
and $\xi_{2}=\xi: AffTA\longrightarrow Aff\alpha_{0}[\textbf{1}_{A}]B\alpha_{0}[\textbf{1}_{A}]$ (note that since $AffTQBQ$ only depends on the unitary equivalence class of $Q$,  it is convenient to denote it as $AffT[Q]B[Q]$), then we will denote $\xi_{1}$ by  $\xi|_{([p],\alpha[p])}$. Even for the general case, we can also write $\xi_{1}=\xi_{2}|_{([p],\alpha[p])}$, when $p<q$ as above.

\noindent\textbf{2.41.}~~As in 2.40, let $A=\oplus_{i=1}^{n}A^{i}$, $B=\oplus_{j=1}^{m}B^{j}$ and $p<q\in A$, $\overline{p}<\overline{q}\in B$, with $\alpha_{0}[p]=[\overline{p}]$ and $\alpha_0[q]=[\overline{q}]$. If $$\gamma_{1}: U(pAp)/\widetilde{SU}(pAp)\longrightarrow U(\overline{p}B\overline{p})/\widetilde{SU}(\overline{p}B\overline{p})$$ is compatible with
$$\gamma_{2}: U(qAq)/\widetilde{SU}(qAq)\longrightarrow U(\overline{q}B\overline{q})/\widetilde{SU}(\overline{q}B\overline{q}),$$ then $\gamma_{1}$ is completely determined by $\gamma_{2}$ (since both maps $$U(pAp)/\widetilde{SU}(pAp)\longrightarrow U(qAq)/\widetilde{SU}(qAq),~~~~ U(\overline{p}B\overline{p})/\widetilde{SU}(\overline{p}B\overline{p})\longrightarrow U(\overline{q}B\overline{q})/\widetilde{SU}(\overline{q}B\overline{q})$$
are injective). Therefore we can denote $\gamma_{1}$ by $\gamma_{2}|_{([p],\alpha[p])}$.

\noindent\textbf{2.42.}~~Let us point out that, in 2.40 and 2.41, if $A\in A\mathcal{HD}$ and $B\in A\mathcal{HD}$, $\xi_{1}$ is not completely determined by $\xi_{2}$ and $\gamma_{1}$ is not completely determined by $\gamma_{2}$.

\vspace{0.3in}



\noindent\textbf{\S3. The counter example}

\noindent\textbf{3.1.}~~In this section, we will present an example of $A\T$ algebras to prove that  $Inv'(A)$ or $Inv(A)$ is not completely determined by $Inv^0(A)$. That is, the Hausdorffifized  algebraic $K_1$ group $\{U(pAp)/\overline{DU(pAp)}\}_{p\in proj(A)}$  or $\{U(pAp)/\widetilde{SU}(pAp)\}_{p\in proj(A)}$ with  the corresponding compatibilities are indispensable   as a  part of the invariant for $Inv'(A)$ or  $Inv(A)$. This is one of the essential  differences between the  simple $C^{*}$-algebras and the $C^{*}$-algebras with the ideal
property. In fact, for all the unital $C^{*}$-algebras $A$ satisfy a  reasonable condition (e.g.,   $\rho(\mathcal{P}_{1}K_{0}(A))=\rho K_{0}(A)$ and $\overline{DU_{0}(A)}=\overline{DU(A)}$),  we have $$U(pAp)/\overline{DU(pAp)}\cong AffTpAp/\overline{\rho K_{0}(pAp)}\oplus K_{1}(pAp),~~~\mbox{and}$$
$$ ~~~~
U(pAp)/\widetilde{SU}(pAp)\cong AffTpAp/\widetilde{\rho K_{0}}(pAp)\oplus K_{1}(pAp)/tor\;K_{1}(pAp),$$ i.e., the metric groups $U(pAp)/\overline{DU(pAp)}$ and  $U(pAp)/\widetilde{SU}(pAp)$ themselves  are completely determined by $AffTpAp$ and $K_{1}(pAp)$, which are included in othe parts of the invariants i.e., there are decided by  $Inv^0(A)$,
but  the compatibilities  make the difference. The point is that the above isomorphisms are
 not natural and therefor the isomorphisms corresponding
to the  cutting down algebras $pAp$ and $qAq$ ($p<q$) may not be chosen to be compatible.

As pointed out in 2.37, $Inv'(A)\cong Inv'(B)$ implies $Inv(A)\cong Inv(B)$. For the $C^*$ algebras $A$ and $B$ constructed in this paper, we only need to prove $Inv^0(A)\cong Inv^0(B)$ but $Inv(A)\not\cong Inv(B)$. Consequently, $Inv'(A)\not\cong Inv'(B)$.

\noindent\textbf{3.2.}~~Let $p_{1}=2,~~~p_{2}=3,~~~p_{3}=5,~~~p_{4}=7,~~~p_{5}=11,\cdots,~~~p_{n}$ be the n-th prime number, let $1<k_{1}<k_{2}<k_{3}<\cdots$ be a sequence of positive integers.
Let
$$A_{1}=B_{1}=C(S^{1}),$$
$$ A_{2}=B_{2}=M_{p_{1}^{k_{1}}}(C[0,1])\oplus M_{p_{1}^{k_{1}}}(C(S^{1}))=A_{1}^{1}\oplus A_{1}^{2}=B_{1}^{1}\oplus B_{1}^{2},$$
$$ A_{3}=B_{3}=M_{p_{1}^{k_{1}}p_{1}^{k_{2}}}(C[0,1])\oplus M_{p_{1}^{k_{1}}p_{2}^{k_{2}}}(C[0,1])\oplus M_{p_{1}^{k_{1}}p_{2}^{k_{2}}}(C(S^{1})),$$
$$ A_{4}=B_{4}=M_{p_{1}^{k_{1}}p_{1}^{k_{2}}p_{1}^{k_{3}}}(C[0,1])\oplus M_{p_{1}^{k_{1}}p_{2}^{k_{2}}p_{2}^{k_{3}}}(C[0,1])\oplus M_{p_{1}^{k_{1}}p_{2}^{k_{2}}p_{3}^{k_{3}}}(C[0,1])\oplus M_{p_{1}^{k_{1}}p_{2}^{k_{2}}p_{3}^{k_{3}}}(C(S^{1})).$$
In general, let\\
$A_{n}=B_{n}=\bigoplus\limits_{i=1}\limits^{n-1}M_{p_{1}^{k_{1}}p_{2}^{k_{2}}\cdots p_{i}^{k_{i}}p_{i}^{k_{i+1}}\cdots p_{i}^{k_{n-1}}}(C[0,1])\oplus M_{p_{1}^{k_{1}}p_{2}^{k_{2}}\cdots p_{n-1}^{k_{n-1}}}(C(S^{1}))$

~~~~~~~~$=\bigoplus\limits_{i=1}\limits^{n-1}M_{\prod\limits_{j=1}\limits^{i}p_{j}^{k_{j}}\cdot\prod\limits_{j=i+1}\limits^{n-1}p_{i}^{k_{j}}}
(C[0,1])\oplus M_{\prod\limits_{i=1}\limits^{n-1}p_{i}^{k_{i}}}(C(S^{1})).$\\
For $1\leq i\leq n-1$,  let $[n,i]=\prod\limits_{j=1}\limits^{i}p_{j}^{k_{j}}\cdot\prod\limits_{j=i+1}\limits^{n-1}p_{i}^{k_{j}}$
and $[n,n]=[n,n-1]$. Then $$A_{n}=B_{n}=\bigoplus\limits_{i=1}\limits^{n-1}M_{[n,i]}(C[0,1])\oplus M_{[n,n]}(C(S^{1})).$$ (Note that last two blocks have same size $[n,n]=[n,n-1]$.)

Note that $[n+1,i]=[n,i]\cdot p_{i}^{k_{n}}$ for  all  $i\in\{1,2,\cdots,n-1\}$ and $[n+1,n+1]=[n+1,n]=[n,n]\cdot p_{n}^{k_{n}}$.

\noindent\textbf{3.3.}~~Let $\{t_{n}\}_{n=1}^{\infty}$ be a dense subset of $[0,1]$ and $\{z_{n}\}_{n=1}^{\infty}$ be a dense subset of $S^{1}$.

In this subsection, we will define the connecting homomorphisms $$\phi_{n,n+1}: A_{n}\longrightarrow A_{n+1}~~~\mbox{ and}~~~~~\psi_{n,n+1}: B_{n}\longrightarrow B_{n+1}.$$
For $i\leq n-1$, define
$\phi_{n,n+1}^{i,i}=\psi_{n,n+1}^{i,i}: ~~M_{[n,i]}(C[0,1])\longrightarrow M_{[n+1,i]}(C[0,1])(=M_{[n,i]\cdot p_{i}^{k_{n}}}(C[0,1]))$ by
$$\phi_{n,n+1}^{i,i}(f)(t)=\psi_{n,n+1}^{i,i}(f)(t)=diag(\underbrace{f(t),f(t),\cdots,f(t)}\limits_{p_{i}^{k_{n}}-1},f(t_{n})),~~ \forall
f\in M_{[n,i]}(C[0,1]).$$
Define $\phi_{n,n+1}^{n,n+1}=\psi_{n,n+1}^{n,n+1}:~~ M_{[n,n]}(C(S^{1}))\longrightarrow M_{[n+1,n+1]}(C(S^{1}))=M_{[n,n]\cdot p_{n}^{k_{n}}}(C(S^{1}))$ by
$$\phi_{n,n+1}^{n,n+1}(f)(z)=\psi_{n,n+1}^{n,n+1}(f)(z)=diag(f(z),\underbrace{f(z_{n}),f(z_{n}),\cdots,f(z_{n})}\limits_{p_{n}^{k_{n}}-1}),~~~\forall f\in
M_{[n,n]}(C(S^{1})).$$
But $\phi_{n,n+1}^{n,n}$ and $\psi_{n,n+1}^{n,n}$ are defined differently---this is the only non-equal component of $\phi_{n,n+1}$ and $\psi_{n,n+1}$.

Let $l=p_{n}^{k_{n}}-1$, then
$$\phi_{n,n+1}^{n,n}(f)(t)=diag(f(e^{2\pi it}),f(e^{-2\pi it}),f(e^{2\pi i\frac{1}{l}}),\cdots,f(e^{2\pi i\frac{l-1}{l}}))$$
$$\psi_{n,n+1}^{n,n}(f)(t)=diag(f(e^{2\pi il_{n}t}),f(e^{-2\pi i\frac{0}{l}}),f(e^{2\pi i\frac{1}{l}}),\cdots,f(e^{2\pi i\frac{l-1}{l}}))$$
for any $f\in M_{[n,n]}(C(S^{1}))$, where $l_{n}=4^{n}\cdot[n+1,n]\in\mathbb{N}$.

Let all other parts $\phi_{n,n+1}^{i,j},\psi_{n,n+1}^{i,j}$ of  $\phi_{n,n+1},\psi_{n,n+1}$ (except $i=j\leq n$ or $i=n,j=n+1$, defined above)
be zero.

Note that all $\phi_{n,n+1}^{i,j},\psi_{n,n+1}^{i,j}$ are either injective or zero.

Let $A=lim(A_{n},\phi_{n,m}),B=lim(B_{n},\psi_{n,m})$. Then it follows from the density of the sets $\{t_{n}\}_{n=1}^{\infty}$ and $\{z_{n}\}_{n=1}^{\infty}$ that both $A$ and $B$ have the ideal property (see the characterization theorem for $AH$ algebras with the ideal property [Pa]).

\noindent\textbf{Proposition 3.4.}~~ There is an isomorphism between $Inv^{0}(A)$ and $Inv^{0}(B)$ (see 2.10)---that is, there is  an isomorphism
$$\alpha: (\underline{K}(A),\underline{K}(A)^{+},\Sigma A)\longrightarrow(\underline{K}(B),\underline{K}(B)^{+},\Sigma B),$$ which is compatible with
Bockstein operations, and for pairs $(p,q)$ with $p\in\Sigma A, q\in\Sigma B$ and  $\alpha([p])=[q]$, there  are  associated unital positive linear maps
$$\xi^{p,q}: AffT(pAp)\longrightarrow AffT(qBq)$$ which are compatible in the sense of 2.9 (see diagram (2.A) in 2.9).
\begin{proof}
Since $KK(\phi_{n,m})=KK(\psi_{n,m})$ and $\phi_{n,m}\sim_{h}\psi_{n,m}$, the identity maps $\eta_{n}: A_{n}\longrightarrow B_{n}$
induce  a shape equivalence between $A=lim(A_{n},\phi_{n,m})$ and $B=lim(B_{n},\psi_{n,m})$, and therefore induce an isomorphism
$$\alpha: (\underline{K}(A),\underline{K}(A)^{+},\Sigma A)\longrightarrow(\underline{K}(B),\underline{K}(B)^{+},\Sigma B).$$
Note that $\phi_{n,n+1}^{i,i}=\psi_{n,n+1}^{i,i}$ for $i\leq n-1$, $\phi_{n,n+1}^{n,n+1}=\psi_{n,n+1}^{n,n+1}$, and
$$\|AffT\phi_{n,n+1}^{n,n}(f)-AffT\psi_{n,n+1}^{n,n}(f)\|\leq\frac{2}{p_{n}^{k_{n}}}\|f\|$$ (see the definition of $\phi_{n,n+1}$ and $\psi_{n,n+1}$).
Therefore,
 $$AffT\eta_{n}: AffTA_{n}\longrightarrow AffTB_{n}~~~\mbox{and }
 ~~~AffT\eta_{n}^{-1}: AffTB_{n}\longrightarrow AffTA_{n}$$
 induce the approximately intertwining diagram
$$
\begin{array}{ccccccccccccc}
AffTA_{1}&\lr& AffTA_{2}&\lr& \cdots&\lr&
 AffTA\vspace{.07in}\\
~\Big\downarrow ~\Big\uparrow& &
~\Big\downarrow ~\Big\uparrow&
\hspace{-.15in}
 \hspace{-.15in}
\hspace{-.15in} \\
AffTB_{1}&\lr& AffTB_{2}&\lr& \cdots&\lr&
 AffTB\vspace{.07in}\\
\end{array}
$$
in the sense of Elliott [Ell1]. Therefore,  there is a unital positive isomorphism
$$\xi: AffTA\longrightarrow AffTB.$$
Also,  for any projection $[P]\in K_{0}(A)$, there is a projection $P_{n}\in A_{n}=B_{n}$ (for $n$ large enough) with $P_{n}^{i}=diag(1,\cdots,1,0,\cdots,0)\in M_{[n,i]}(C(X_{n,i}))$, where $X_{n,i}=[0,1]$ for $i\leq n-1$, and $X_{n,n}=S^{1}$, such that
$\phi_{n,\infty}([P_{n}])=[P]\in K_{0}(A)$. Note that for any  constant functions $f\in A_{n}^{i}=B_{n}^{i}$ (e.g., $P_{n}^{i}$ above) and for any $j$,  $\phi^{i,j}_{n,n+1}(f)$ and $\psi^{i,j}_{n,n+1}(f)$ are still constant functions and $\phi^{i,j}_{n,n+1}(f)=\psi^{i,j}_{n,n+1}(f)$.
That is, we have
$$\phi_{n,n+1}(P_{n})=\psi_{n,n+1}(P_{n}) ~~(\mbox{denoted by}~   P_{n+1})~~~~~~~\mbox{and }$$
$$\phi_{n,m}(P_{n})=\psi_{n,m}(P_{n})~~(\mbox{denoted by}~P_{m}).$$
 Let $P_{\infty}=\phi_{n,\infty}(P_{n})$ and $Q_{\infty}=\psi_{n,\infty}(P_{n})$.  Then the identity maps $\{\eta_{m}\}_{m>n}$ also induce the following approximate intertwining diagram:
$$
\begin{array}{ccccccccccccc}
AffT(P_{n}A_{n}P_{n})&\lr& AffT(P_{n+1}A_{n+1}P_{n+1})&\lr& \cdots&\lr&
 AffTP_{\infty}AP_{\infty}\vspace{.07in}\\
~\Big\downarrow ~\Big\uparrow&&
~\Big\downarrow ~\Big\uparrow&
\hspace{-.15in}
 \hspace{-.15in}
\hspace{-.15in} \\
AffT(P_{n}B_{n}P_{n})&\lr& AffT(P_{n+1}B_{n+1}P_{n+1})&\lr& \cdots&\lr&
 AffTQ_{\infty}BQ_{\infty}~,\vspace{.07in}

\end{array}
$$
and hence induce a positive linear isomorphism
$$\xi^{[P],\alpha[P]}: AffTP_{\infty}AP_{\infty}\rightarrow AffTQ_{\infty}BQ_{\infty}.$$
(Note that $[P_{\infty}]=[P],\;[Q_{\infty}]=\alpha [P]\;in\;K_{0}(A)\;and\;K_{0}(B)$,\;respectively.)
Evidently those maps are compatible since, they are induced by the same sequence of homomorphisms $\{\eta_{n}\}$ and $\{\eta_{n}^{-1}\}$.\\
\end{proof}

The following Definition 3.5 and Proposition 3.6 are  inspired by [Ell3].

\noindent\textbf{Definition 3.5.}~~ Let $C\!=\!\lim(C_{n},\phi_{n,m})$\;be\;an\;$A\mathcal{HD}$ inductive limit. We say the system $(C_{n},\phi_{n,m})$ has the uniformly varied determinant if for any $C^{i}_{n}=M_{[n,i]}(C(S^{1}))$\;(that is,  $C^{i}_{n}$\;has spectrum\;$S^{1}$)\;and $C^{j}_{n+1}$ and $f\in C_n^i$ defined by
$$f(z)=\left(
        \begin{array}{cccc}
          z &  &  &  \\
           & 1 &  &  \\
           &  & \ddots &  \\
           &  &  & 1 \\
        \end{array}
      \right)_{[n,i]\times[n,i]} ~~~~~\forall z\in S^1,$$
we have that det($\phi^{i,j}_{n,n+1}(f)(x))=$ constant for $x\in Sp(C^{j}_{n+1})\neq S^{1}$ or det($\phi^{i,j}_{n,n+1}(f)(z))=\lambda z^{k}$ ($\lambda \in\mathbb{C}$) for $z\in Sp(C^{j}_{n+1})=S^{1}$, where $j$ satisfy $\phi^{i,j}_{n,n+1}\not=0$ and the  determinant is taken inside $\phi^{i,j}_{n,n+1}(1_{C^{i}_{n}})C^{j}_{n+1}\phi^{i,j}_{n,n+1}(1_{C^{i}_{n}})$.

\noindent\textbf{Proposition 3.6.}~~ If the inductive limit system $C=(C_{n},\phi_{n,m})$\;has the uniformly varied determinant, then for any elements  $[p]\in \sum C$, there are a splitting maps
$$K_{1}(pCp)/tor\;K_{1}(pCp)\xrightarrow{S_{pCp}}U(pCp)/\widetilde{SU}(pCp)$$ of the exact sequences
$$0\rightarrow AffT{pCp}/\widetilde{\rho K_{0}}(pCp)\rightarrow U(pCp)/\widetilde{SU}(pCp)\xrightarrow{\pi_{pCp}}K_{1}(pCp)/tor\;K_{1}(pCp)\rightarrow0$$
(that is, $\pi_{pCp}\circ S_{pCp}=id\;on\;K_{1}(pCp)/tor\;K_{1}(pCp)$) such that the system of maps $\{S_{pCp}\}_{[p]\in \sum C}$ are compatible in the following sense: if $p<q$, then the following diagram commutes
$$(3.6.A)~~~~~~~~~~~~~~~~~~~~~~~~~~~\CD
 K_{1}(pCp)/tor\;K_{1}(pCp) @>S_{pCp}>> U(pCp)/\widetilde{SU}(pCp) \\
  @V  VV @V  VV  \\
 K_{1}(qCq)/tor\;K_{1}(qCq) @>S_{qCq}>> U(qCq)/\widetilde{SU}(qCq),
\endCD~~~~~~~~~~~$$
where the vertical maps are  induced by the inclusions $pCp\longrightarrow qCq$.
\begin{proof}
 Fix $p\in C$. Let $x\in K_{1}(pCp)/tor\;K_{1}(pCp)$. There  exist a $C_{n}$ and $p_{n}\in C_{n}$ such that $ [\phi_{n,\infty}(p_{n})]=[p]\in K_{0}(C)$. Without lose of generality, we can assume $\phi_{n,\infty}(p_{n})=p$. By increasing $n$ if necessary, we can assume that there is an element  $x_{n}\in K_{1}(p_{n}C_np_{n})/tor\;K_{1}(p_{n}C_np_{n})$, such that $(\phi_{n,\infty})_*(x_{n})=x \in K_{1}(pCp)/tor\;K_{1}(pCp)$.

Write $p_nC_np_n= D=\oplus D^i$. Let $I=\{i~|~~Sp(D^{i})=S^{1}\}$. For $i\in I$,  $D^i$ can be identified with $M_{l_i}(C(S^1))$.
Let $u_i\in D^i$ be defined by
$$u_{i}(z)=\left(
        \begin{array}{cccc}
          z &  &  &  \\
           & 1 &  &  \\
           &  & \ddots &  \\
           &  &  & 1 \\
        \end{array}
      \right)_{l_i\times l_i}~~~~\forall z\in S^1, $$
    which represents  the standard generator of $K_{1}(D^i).$
Then $x_{n}$ can be represented by $$u=\mathop{\bigoplus}\limits_{i\in I}u_{i}^{k_{i}}\oplus\mathop{\bigoplus}\limits_{j\notin I}\textbf{1}_{D^{j}}~~~\in~~~\mathop{\bigoplus}\limits_{i\in I}D^{i}\oplus\mathop{\bigoplus}\limits_{j\notin I}D^{j}=D\subseteq p_nC_{n}p_n.$$ Define $S(x)=[\phi_{n,\infty}(u)]\in U(pCp)/\widetilde{SU}(pCp)$.
Note that all unitaries with constant determinants  are  in $\widetilde{SU}$, and that the inductive system  has the uniformly  varied determinant, it is routine to verify that $S(x)$ is well defined and the system $\{S_{pCp}\}_{[p]\in \sum C}$ makes the diagram (3.6.A) commute. \\
\end{proof}


\noindent\textbf{3.7.}~~Let ${\cal A}$ be a unital $C^*$-algebra. Then $AffT{\cal A}$ is a real Banach space with quotient space
$AffT{\cal A}/\widetilde{\rho K_0}({\cal A})$. Let us use $\|\cdot\|^{\sim}$ to denote the quotient norm. Note that $\widetilde{\lambda}_{\cal A}$ identifies $U_{tor}({\cal A})/\widetilde{SU}({\cal A})$ with $AffT{\cal A}/\widetilde{\rho K_0}({\cal A})$. In this way, $U_{tor}({\cal A})/\widetilde{SU}({\cal A})$ is regarded as a real Banach space, whose norm is also denoted by $\|\cdot\|^{\sim}$. In general, we have
$$U({\cal A})/\widetilde{SU}({\cal A})\cong U_{tor}({\cal A})/\widetilde{SU}({\cal A})~\times~K_1({\cal A})/torK_1({\cal A});$$
but the identification is not canonical. Even though $U({\cal A})/\widetilde{SU}({\cal A})$ is not a Banach  space, it is an Abelian group: for $[u], [v] \in U({\cal A})/\widetilde{SU}({\cal A})$, define $[u]-[v]=[uv^*]$.

The  norm $\|\cdot\|^{\sim}$ is related
to the metrices $\widetilde{d}_{\cal A}$ (on $AffT{\cal A}/\widetilde{\rho K_0}({\cal A})$; see 2.26) and $\widetilde{D}_{\cal A}$ (on $U_{tor}({\cal A})/\widetilde{SU}({\cal A})$; see 2.30) as below.  Let $\varepsilon <1$.  For any $f, g\in AffT{\cal A}/\widetilde{\rho K_0}({\cal A})$,
$$\|f-g \|^{\sim}<\frac{\varepsilon}{2\pi}~~ \Longrightarrow ~~\widetilde{d}_{\cal A}(f,g)<\varepsilon~~ \Longrightarrow~~
\|f-g \|^{\sim}<\frac{\varepsilon}{4}. $$
And for any $[u], [v]\in U({\cal A})/\widetilde{SU}({\cal A})$ with $[u]-[v]=[uv^*]\in U_{tor}({\cal A})/\widetilde{SU}({\cal A})$,
$$\|[u]-[v]\|^{\sim}<\frac{\varepsilon}{2\pi}~~ \Longrightarrow ~~\widetilde{D}_{\cal A}([u], [v])<\varepsilon~~ \Longrightarrow~~
\|[u]-[v] \|^{\sim}<\frac{\varepsilon}{4}. $$

For ${\cal A}=PM_{l}(C(X))P\in\mathcal{HD}$ or  ${\cal A}=M_{l}(I_{k})$ (at this case we also denote $[0,1]$ by $X$), there are canonical identification (see 2.39)
$$U_{tor}({\cal A})/\widetilde{SU}({\cal A})\cong AffT{\cal A}/\widetilde{\rho K_0}({\cal A}) \cong C(X,\R)/\{constant\; functions\}. $$
Choose a base point $x_0\in X$. Let $C_{x_0}(X, \R)$ be the set of functions $f\in C(X, \R)$ with $f(x_0)=0$. Then $ C(X,\R)/\{constant\; functions\}\cong C_{x_0}(X, \R)$.
For $[f]\in AffT{\cal A}/\widetilde{\rho K_0}({\cal A})$ (or  $[f]\in U_{tor}({\cal A})/\widetilde{SU}({\cal A})$) identified with a function $f\in C_{x_0}(X,\R)$, we have
$$\|[f]\|^{\sim}= \frac{1}{2}\big(\max_{x\in X} (f(x))-\min_{x\in X} (f(x))\big),$$
(rather than $sup_{x\in X}\{|f(x)|\}$).

In the above case, if $p\in {\cal A}$ is a non zero projection,  then $U_{tor}(p{\cal A}p)/\widetilde{SU}(p{\cal A}p)\cong AffT(p{\cal A}p)/\widetilde{\rho K_0}(p{\cal A}p)$  is also identified with $C_{x_0}(X,\R)$.  Consider the inclusion map $\imath: pAp \to A$.  Then the map $\imath_*$ as map from $U_{tor}(p{\cal A}p)/\widetilde{SU}(p{\cal A}p)\cong AffT(p{\cal A}p)/\widetilde{\rho K_0}(p{\cal A}p)$ to $U_{tor}({\cal A})/\widetilde{SU}({\cal A})$
 can be described as below: if \\ $u\in U_{tor}(p{\cal A}p)/\widetilde{SU}(p{\cal A}p)\cong AffT(p{\cal A}p)/\widetilde{\rho K_0}(p{\cal A}p)$ is identified with $f\in C_{x_0}(X,\R)$, then \\ $\imath_*(u)\in U_{tor}({\cal A})/\widetilde{SU}({\cal A})$ is identified with  $\frac{rank(p)}{rank(\one_{cal A})}f$. But $\imath^{\natural}$ is the identity map from $U_{tor}(p{\cal A}p)/\widetilde{SU}(p{\cal A}p)\cong AffT(p{\cal A}p)/\widetilde{\rho K_0}(p{\cal A}p)$ to itself (not to $U_{tor}({\cal A})/\widetilde{SU}({\cal A})$).

\noindent\textbf{3.8.}~~ It is easy to see that $K_1(A)=K_1(B)=\Z$.

In the definition of $A_n=\oplus_{i=1}^nA_n^i$, only one block $A_n^n=M_{[n,n]}(C(S^1))$ has spectrum $S^1$, and only two partial maps  $\phi^{n,j}_{n,n+1}$ for $ ~j=n,~ j=n+1$  (of $\phi_{n,n+1}$ from $A_n^n$) are nonzero.  Let $f\in A_n^n$ be defined as in Definition 3.5. Then det$(\phi^{n,n+1}_{n,n+1}(f)(z))=z$ and det$(\phi^{n,n}_{n,n+1}(f)(t))=e^{2\pi it}e^{-2\pi it}e^{2\pi i\frac1l}e^{2\pi i\frac2l}\cdots e^{2\pi i\frac{l-1}l}=\pm 1$ (see 3.3).
So the inductive limit system $(A_{n},\phi_{n,m})$ has the uniformly varied determinant, and therefore the limit algebra $A$ has compatible splitting maps
$S_{p}: K_{1}(pAp)\rightarrow U(pAp)/\widetilde{SU}(pAp)$.

We will prove that $B=lim(B_{n},\psi_{n,m})$ does not have such compatible system of splitting  maps  $\{K_{1}(pBp)\longrightarrow U(pBp)/\widetilde{SU}(pBp)\}_{[p]\in\sum B}$.

Before  proving  the above fact, let us describe the $K_0$-group of $A$ and $B$.
Let $$G_{1}=\{\frac{m}{p_{1}^{l}}~|~~m\in\mathbb{Z},l\in\mathbb{Z_{+}}\},$$
 $$\;\;\;\;\;G_{2}=\{\frac{m}{p_{1}^{k_{1}}p_{2}^{l}}~|~~m\in\mathbb{Z},l\in\mathbb{Z_{+}}\},$$
 $$\;\;\;\;\;\;\;\;\;\;G_{3}=\{\frac{m}{p_{1}^{k_{1}}p_{2}^{k_{2}}p_{3}^{l}}~|~~m\in\mathbb{Z},l\in\mathbb{Z_{+}}\},$$
$$\vdots$$
$$\;\;\;\;\;\;\;\;\;\;\;\;\;\;\;\;\;\;\;\;G_{n}=\{\frac{m}{p_{1}^{k_{1}}p_{2}^{k_{2}}\ldots p_{n-1}^{k_{n-1}}p_{n}^{l}}~|~~m\in\mathbb{Z},l\in\mathbb{Z_{+}}\},$$
$$\;\;\;\;\;\;\;\;\;\;\;\;\;\;\;G_{\infty}=\{\frac{m}{p^{k_{1}}_{1},p^{k_{2}}_{2},\cdots,p^{k_{t}}_{t}}~|~~t\in \mathbb{Z_{+}}, m\in\mathbb{Z}\},$$
where $p_{1}=2,p_{2}=3, \cdots, p_i,\cdots $ and $k_1,k_2,\cdots, k_i\cdots $ are defined in 3.2. Then
 $$K_{0}(A)=K_{0}(B)=\{(a_{1},a_{2},\cdots,a_{n},\cdots)\in \prod\limits_{n=1}\limits^{\infty}G_{n}~|~~\exists N~~\mbox{such that}~~a_{N}=a_{N+1}=\cdots\in \mathbb{Q}\}\triangleq{\tilde G}.$$
  Furthermore, their positive cones consist  of the elements whose
  coordinates are non-negative,
   and their order units are $[\one_A]=[\one_B]=(1,1,\cdots,1,\cdots)\in \prod\limits_{n=1}\limits^{\infty}G_{n}$.
    Let $$\alpha_0: (K_0(A), K_0(A)^+,[\one_A])=\big({\tilde G},{\tilde G}^+,(1,1,\cdots,1,\cdots)\big)  \to K_0(B), K_0(B)^+,[\one_B])=\big({\tilde G},{\tilde G}^+,(1,1,\cdots,1,\cdots)\big)$$
    be a scaled ordered isomorphism. Then  $\alpha_0 \big((1,1,\cdots,1,\cdots)\big)=(1,1,\cdots,1,\cdots)$. Note that an element  $x\in {\tilde G}$ is divisible by   power $p_1^n$ (for any $n$) of the first prime number $p_1=2$ if and only if $x=(t,0,0,\cdots, 0,\cdots)\in G_1\subset {\tilde G}$. Hence $\alpha_0\big((1,0,0,\cdots, 0,\cdots)\big)=(t,0,0,\cdots, 0,\cdots)$ for some $t\in G_1$ with $t>0$. Hence
  $$\alpha_0\big(0,1,1,\cdots, 1,\cdots)\big)=(1-t,1,1,\cdots,1,\cdots).$$
 Since $\alpha_0$ preserves the positive cone, we have  $1-t\geq 0$ which implies $t\leq 1$. On the other hand, $(\alpha_0)^{-1}$ takes $(1,0,0,\cdots,0,\cdots)$ to $(1/t,0,0,\cdots,0,\cdots)$. But $(\alpha_0)^{-1}$ also preserves the positive cone. Symmetrically, we get $t\geq 1$. That is, $\alpha_0\big((1,0,0,\cdots, 0,\cdots)\big)=(1,0,0,\cdots, 0,\cdots)$. Similarly, using  the fact that $G_k$ is the subgroup of all elements in ${\tilde G}$ which can be divisible by any power of $p_k$---the $k^{th}$ prime number, we can prove that
 $$\alpha_0\big((\underbrace{0,\cdots,0}_{k-1},1,0,\cdots,0,\cdots)\big)=
 (\underbrace{0,\cdots,0}_{k-1},1,0,\cdots,0,\cdots)\in G_k\subset{\tilde G}.$$
 That is,  $\alpha_0$ is the identity on ${\tilde G}$.

Note that $Sp(A)=Sp(B)$ is the  one point compactification of $\{1,2,3\cdots\}$---or,  in other words,
$\{1,2,3\cdots,\infty\}$. If we let $I_{n}$ (or $J_{n}$) be the primitive ideal $A$ (or $B$) corresponding to $n$ (including $n=\infty$), then
$$K_{0}(A/I_{n})=K_{0}(B/J_{n})=G_{n}.$$

Note also that if $m'>m>n \in \N$, then $\phi_{m,m'}(A_m^n)\subset A_{m'}^n$ and $\psi_{m,m'}(B_m^n)\subset B_{m'}^n$. Hence \\
$A/I_n=\lim_{n<m\to \infty} (A_m^n, \phi_{m,m'}|_{A_m^n})$ (and $B/J_n=\lim_{n<m\to \infty} (B_m^n, \psi_{m,m'}|_{B_m^n})$ resp.) are   ideals of $A$ (and $B$ resp.). But $A/I_{\infty}$ (or
$B/J_{\infty}$) is not an ideal of $A$ (or $B$).

  Let $\alpha: (\underline{K}(A),\underline{K}(A)^{+},\Sigma A)\longrightarrow(\underline{K}(B),\underline{K}(B)^{+},\Sigma B)$
 be an isomorphism.
By 3.8 the induced map $\alpha_0$ on $K_0$ group is identity, when both
$K_0(A)$ and $K_0(B)$ are identified with ${\tilde G}$ as scaled ordered groups. That is, $\alpha_0$ is the same as the $\alpha_0$  induced by the shape equivalence in  the proof of Proposition 3.4. In particular, if
there is an isomorphism
 $\wedge: A\longrightarrow B$, then for all  $i\leq n-1$,
$\wedge_*[(\phi_{n,\infty}(\textbf{1}_{A_{n}^{i}}))]=
[\psi_{n,\infty}(\textbf{1}_{B_{n}^{i}})]$. This implies $\wedge(\phi_{n,\infty}(\textbf{1}_{A_{n}^{i}}))=
\psi_{n,\infty}(\textbf{1}_{B_{n}^{i}})$, since $\psi_{n,\infty}(\textbf{1}_{B_{n}^{i}})=\one_{B/I_i}$, which is in the center of $B$ (any element in the center of the $C^*$-algebra can only unitary equivalent to itself).
Hence it is also true that $\wedge(\phi_{n,\infty}(\textbf{1}_{A_{n}^{i}}))=
\psi_{n,\infty}(\textbf{1}_{B_{n}^{i}})$ for  $i=n$.

\noindent\textbf{3.9.}~~Let $P_{1}=1_{B}=\psi_{1,\infty}(\textbf{1}_{B_{1}}), ~~P_{2}=\psi_{2,\infty}(\textbf{1}_{B_{2}^{2}}), ~~P_{3}=\psi_{3,\infty}(\textbf{1}_{B_{3}^{3}}),
\cdots,~~P_{n}=\psi_{n,\infty}(\textbf{1}_{B_{n}^{n}}),\cdots.$ Then \\ $P_{1}>P_{2}>\cdots>P_{n}\cdots.$
We will prove that there are no
splitings $$K_{1}(P_{n}BP_{n})\longrightarrow U(P_{n}BP_{n})/\widetilde{SU}(P_{n}BP_{n})$$
which are compatible for  all pairs of projections $P_n>P_m$ (see diagram (3.6.A)), in the
next subsection. Before doing so,  we need some preparations.

Set $Q_{1}=P_{1}-P_{2},~~Q_{2}=P_{2}-P_{3},\cdots,~~Q_{n}=P_{n}-P_{n+1}.$ Then for each $n$, we have the inductive limit
$$Q_{n}BQ_{n}=\lim\limits_{m\rightarrow\infty}(B_{m}^{n},\psi_{m,m'}^{n,n}),$$
(note that for $m>n$, $\psi_{m,m+1}^{n,j}=0$ if $j\neq n$), which is the quotient algebra corresponding to the primitive ideal of $n\in Sp(B)=\{1,2,3\cdots,\infty\}$. Note that $Q_{n}BQ_{n}$ is a simple $AI$ algebra. The inductive limit of the $C^*$-algebras
$$B_{n+1}^{n}\longrightarrow B_{n+2}^{n}\longrightarrow B_{n+3}^{n}\longrightarrow\cdots\longrightarrow Q_{n}BQ_{n}$$
induces the inductive limit of the ordered Banach spaces
$$AffTB_{n+1}^{n}\xrightarrow{\xi_{n+1,n+2}}AffTB_{n+2}^{n}\xrightarrow{\xi_{n+2,n+3}}\cdots\rightarrow AffTQ_{n}BQ_{n},$$
whose connecting maps $\xi_{m,m+1}:  C_{\mathbb{R}}([0,1])\longrightarrow C_{\mathbb{R}}([0,1])$ (for $m>n$) satisfy that
$$\|\xi_{m,m+1}(f)-f\|\leq\frac{1}{p_{n}^{k_{m}}}\|f\|,~~~~~\forall f\in C_{\mathbb{R}}[0,1],~~ m>n.$$ Hence we have the following approximate intertwining diagram
$$
\begin{array}{ccccccccccccc}
C_{\mathbb{R}}[0,1]&\xrightarrow{\xi_{n,n+1}}& C_{\mathbb{R}}[0,1]&\xrightarrow{\xi_{n+1,n+2}}& C_{\mathbb{R}}[0,1]&\lr&
 \cd&\lr& AffTQ_{n}BQ_{n}\vspace{.07in}\\
~\Big\downarrow ~\Big\uparrow&&
~\Big\downarrow ~\Big\uparrow&&~\Big\downarrow ~\Big\uparrow&
\hspace{-.15in}
 \hspace{-.15in}
\hspace{-.15in} \\
C_{\mathbb{R}}[0,1]&\xrightarrow{\;\;\;\;id\;\;\;\;}& C_{\mathbb{R}}[0,1]
&\xrightarrow{\;\;\;\;id\;\;\;\;}& C_{\mathbb{R}}[0,1]&\lr&\cd &\lr&C_{\mathbb{R}}[0,1]~~~.\vspace{.17in}\\

\end{array}
$$
Consequently, $AffTQ_{n}BQ_{n}\cong C_{\mathbb{R}}[0,1]$, and the maps $$\xi_{m,\infty}: AffTB_{m}^{n}=C_{\mathbb{R}}[0,1]\longrightarrow AffTQ_{n}BQ_{n}\cong C_{\mathbb{R}}[0,1]$$ (under the identification) satisfy
 $$\|\xi_{m,\infty}(f)-f\|\leq(\frac{1}{p_{n}^{k_{m}}}+\frac{1}{p_{n}^{k_{m+1}}}+\cdots)\|f\|\leq
\frac{1}{4}\|f\|,~~~~~\forall f\in C_{\mathbb{R}}[0,1].$$
Therefore  $\|\xi_{m,\infty}(f)\|\geq\frac{3}{4}\|f\|.$ \\
Note that $\widetilde{\rho K_{0}}(Q_{n}BQ_{n})=\mathbb{R}=\widetilde{\rho K_{0}}(B_{m}^{n})$ consists of constant functions on $[0,1]$. Let
$h\in C_{\mathbb{R}}[0,1]=
AffT(B_{m}^{n}).$
Considering  the element $\xi_{m,\infty}(h)$ as in $AffT(Q_{n}BQ_{n})/\widetilde{\rho K_{0}}(Q_{n}BQ_{n})$, we have $$\|\xi_{m,\infty}(h)\|^{\sim}\geq\frac{1}{2}\cdot\frac{3}{4}(\max\limits_{t\in[0,1]}h(t)-\min\limits_{t\in[0,1]}h(t)),$$
where $\|\cdot\|^{\sim}$ is defined in 3.7.

\noindent\textbf{3.10.}~~We now prove that  no compatible splittings $$S_{n}: K_{1}(P_{n}BP_{n})\longrightarrow U(P_{n}BP_{n})/\widetilde{SU}(P_{n}BP_{n})$$
exists. Suppose such splittings exist. Then consider the generator
$x\in K_{1}(B)=\mathbb{Z}.$

Note that $x\in K_{1}(P_{n}BP_{n})\cong K_{1}(B)$,  for all $P_{n}$. Note also that
the diagram $$\CD
  K_{1}(P_{n+1}BP_{n+1}) @>S_{n+1}>> U(P_{n+1}BP_{n+1})/\widetilde{SU}(P_{n+1}BP_{n+1})\\
  @V id VV @V   \imath_{*} VV \\
   K_{1}(P_{1}BP_{1})@>S_{1}>> U(P_{1}BP_{1})/\widetilde{SU}(P_{1}BP_{1})
\endCD$$
\\
commutes $(P_{1}BP_{1}=B)$. The composition
$$U(P_{n+1}BP_{n+1})/\widetilde{SU}(P_{n+1}BP_{n+1})\mathop{\longrightarrow}\limits^{\imath_*} U(P_{1}BP_{1})/\widetilde{SU}(P_{1}BP_{1})\longrightarrow\bigoplus\limits_{i=1}
\limits^{n}U(Q_{i}BQ_{i})/\widetilde{SU}(Q_{i}BQ_{i})$$
is the zero map. (Note that $Q_{i}BQ_{i}$ is an ideal of $B$ and  is also the   quotient $B/J_{i}$.)
 Consequently, we have
$$(*)~~~~~~~~~~~~~~~~~~~~~~~~~~~~~~~~~~~~~~~~~~~~\pi_{n}^{\natural}(S_1(x))=\pi_{n}^{\natural}( \imath_*S_{n+1}(x))=0,~~~~~~~~~~~~~~~~~~~~~~~~$$
where $\pi_{n}:B\rightarrow Q_{n}BQ_{n}$ is the quotient map.
Let $S_{1}(x)$ be represented by a unitary $u\in U(B)$. Then there are an  $n$ (large enough) and $[u_n]\in U(B_{n})/\widetilde{SU}(B_{n})$, represented by unitary
$u_{n}\in B_{n}$, such that $\psi_{n,\infty}^{\natural}([u_{n}])-S_{1}(x)\in U_{tor}(B_n)/\widetilde{SU}(B_n)$ and
$$\|\psi_{n,\infty}^{\natural}([u_{n}])-S_{1}(x)\|^{\sim}<\frac{1}{16}.$$
Note that $$(\psi_{n,m})_{*}: K_{1}(B_{n})\longrightarrow K_{1}(B_{m})$$
is the identify map from $\mathbb{Z}$ to $\mathbb{Z}$. Let $g\in M_{[n,n]}(C(S^{1}))=B_{n}^{n}$ be defined by
$$g(z)=\left(
    \begin{array}{ccccc}
      z &  &  &  &  \\
       & 1 &  &  &  \\
       &  & 1 &  &  \\
       &  &  & \ddots &  \\
       &  & &  & 1 \\
    \end{array}
  \right)_{[n,n]\times[n,n]~~.}$$
Then $[g^{-1}u_{n}]=0$ in $K_{1}(B_{n})$. By the exactness of the  sequence
$$0\longrightarrow AffTB_{n}/\widetilde{\rho K_{0}}(B_{n})\longrightarrow U(B_{n})/\widetilde{SU}(B_{n})\longrightarrow K_{1}(B_{1})\longrightarrow0,$$ there is an $h\in\bigoplus\limits_{i=1}\limits^{n}C_{\mathbb{R}}[0,1]\oplus C_{\mathbb{R}}(S^{1})=AffTB_{n}$ such that
$$[u_{n}]=[g]\cdot(e^{2\pi ih}\cdot1_{B_{n}})\in U(B_{n})/\widetilde{SU}(B_{n}).$$ Let $\|h\|=M$. Choose $m>n$ such that $4^{m-1}>8M+8$.

Consider
$$\psi_{n,m}^{n,m-1}: B_{n}^{n}=M_{[n,n]}(C(S^{1}))\longrightarrow B_{m}^{m-1}=M_{[m,m-1]}(C([0,1]))$$ which is the composition
$$\psi_{m-1,m}^{m-1,m-1}\circ\psi_{n,m-1}^{n,m-1}: M_{[n,n]}(C(S^{1}))\longrightarrow M_{[m-1,m-1]}(C(S^{1}))\longrightarrow M_{[m,m-1]}(C([0,1])).$$
Let $g'=\psi_{n,m}^{n,m-1}(g)$. We know that
$$g'(t)=\psi_{n,m}^{n,m-1}(g)(t)=\left(
                                               \begin{array}{ccccc}
                                                 e^{2\pi il_{m-1}t} &  &  &  &  \\
                                                  & * &  &  &  \\
                                                  &  & * &  &  \\
                                                  &  &  & \ddots &  \\
                                                  &  &  &  & * \\
                                               \end{array}
                                             \right)_{[m,m-1]\times[m,m-1]},$$
where the $*'s$ represent constant functions on $[0,1]$, and therefore
$$g'=e^{2\pi ih'}~~~~~(mod~~\widetilde{SU}(B_{m}^{m-1})))$$
with $h'(t)=\frac{l_{m-1}}{[m,m-1]}\cdot t\cdot 1_{[m,m-1]}.$ When we identify $U(B_{m}^{m-1})/\widetilde{SU}(B_{m}^{m-1})$ with
$$AffTB_{m}^{m-1}/\widetilde{\rho K_{0}}(B_{m}^{m-1})=C_{\mathbb{R}}[0,1]/\{constants\},$$
$g'$ is identified with $\widetilde{h}\in C_{\mathbb{R}}[0,1]$ with
$$\widetilde{h}(t)=\frac{l_{m-1}}{[m,m-1]}t.$$
Since $\frac{l_{m-1}}{[m,m-1]}\geq8M+8,$ we have
$$\|\widetilde{h}\|^{\sim}=\frac{1}{2}\big(\max\limits_{t\in[0,1]}\widetilde{h}(t)-\min\limits_{t\in[0,1]}\widetilde{h}(t)\big)\geq 4M+4$$
(see 3.7). On the other hand, $$[u_{n}]=[g]+\widetilde{\lambda}_{B_{n}}([h])\in U(B_{n})/\widetilde{SU}(B_{n}),$$ where $[h]\in AffTB_{n}/\widetilde{\rho K_{0}}(B_{n})$ is the
element defined by $h$, and $$\widetilde{\lambda}_{B_{n}}: AffTB_{n}/\widetilde{\rho K_{0}}(B_{n})\longrightarrow U(B_{n})/\widetilde{SU}(B_{n})$$
is the map defined in 2.30 (also see 2.26). Consequently,
$$(\psi_{n,m}^{n,m-1})^{\natural}(u)=AffT\psi_{n,m}^{n,m-1}(h)+\widetilde{h}\triangleq\widetilde{\widetilde{h}}\in AffTB_{m}^{m-1}/\widetilde{\rho K_{0}}(B_{m}^{m-1})\cong U(B_{m}^{m-1})/\widetilde{SU}(B_{m}^{m-1})$$ with
$$\|\widetilde{\widetilde{h}}\|^{\sim}=\frac{1}{2}\big(\max\limits_{t\in[0,1]}\widetilde{\widetilde{h}}(t)-\min\limits_{t\in[0,1]}\widetilde{\widetilde{h}}(t)\big)\geq4,$$
since $\|h\|\leq M.$ Therefore,
$$(\pi_{m-1}\circ\psi_{n,\infty})^{\natural}(u)\in U(Q_{m-1}BQ_{m-1})/\widetilde{SU}(Q_{m-1}BQ_{m-1})\cong AffT(Q_{m-1}BQ_{m-1})/\widetilde{\rho K_{0}}(Q_{m-1}BQ_{m-1}),$$ satisfies
$$\|(\pi_{m-1}\circ\psi_{n,\infty})^{\natural}(u)\|^{\sim}=\frac{1}{2}\big(\max\limits_{t\in[0,1]}(\pi_{m-1}\circ\psi_{n,\infty})^{\natural}(u)(t)-\min\limits_{t\in[0,1]}(\pi_{m-1}\circ\psi_{n,\infty})^{\natural}(u)(t)\big)\geq\frac{3}{4}\cdot 4=3,$$
where $\pi_{m-1}: B\longrightarrow Q_{m-1}BQ_{m-1}$ is the quotient map.
On the other hand,  $$\pi_{m-1}^{\natural}(S_{1}(x))=0$$ as calculated in ($*$). Recall that
$$\|(\psi_{n,\infty})^{\natural}(u)-S_{1}(x)\|^{\sim}<\frac{1}{16}.$$
We get
$$\|(\pi_{m-1}\circ\psi_{n,\infty})^{\natural}(u)\|^{\sim}<\frac{1}{16}$$ which is a contradiction. This contradiction proves that such system of splittings do not exist. Hence $Inv(A)\ncong Inv(B)$ and $A\ncong B$.

\noindent\textbf{3.11.} One can easily verify that $$AffTA=AffTB=\{(f_{1},f_{2},\cdots,f_{n}\cdots)\in \prod\limits_{n=1}\limits^{\infty}C_{\R}[0,1]~|~~\exists~ r\in  \R~~\mbox{such that}~~f_n(x)~\mbox{conveges to} ~r~\mbox{uniformly}\}.$$
$$\overline{\rho K_{0}(A)}(=\overline{\rho K_{0}(B)})=\{(r_{1},r_{2},\cdots,r_{n},\cdots)\in \prod\limits_{n=1}\limits^{\infty}\R~|~~\exists~ r\in  \R~~\mbox{such that}~~r_n~\mbox{conveges to} ~r~\}\subset AffTA (=AffTB).$$
Since $\overline{\rho K_{0}(A)}(=\overline{\rho K_{0}(B)})$ is already a vector space, we have $\widetilde{\rho K_0}(A)=\overline{\rho K_{0}(A)}$ and $\widetilde{\rho K_0}(B)=\overline{\rho K_{0}(B)}$.
 Therefore
$$U_{tor}(A)/\widetilde{SU}(A)\cong AffTA/\widetilde{\rho K_0}(A)=AffTA/\overline{\rho K_{0}(A)}\cong U_0(A)/\overline{DU(A)}.$$
On the other hand, $U_{tor}(A)=U_0(A)$. Hence $\widetilde{SU}(A)=\overline{DU(A)}$.
Furthermore  the map ${\lambda}_{A}: AffTA/\overline{\rho K_{0}(A)}\longrightarrow U(A)/\overline{DU(A)}$ can be identified with the map
$\widetilde{\lambda}_{A}: AffTA/\widetilde{\rho K_{0}}(A)\longrightarrow U(A)/\widetilde{SU(A)}.$
That is $Inv'(A)=Inv(A)$. Similarly, $Inv(B)=Inv'(B)$.

\noindent\textbf{3.12.} A routine calculation shows (we omit the details) that  for any finite subset $F\subset A_n$, and $\varepsilon>0$, there is an $m>n$ and two finite dimensional  unital sub $C^*$-algebras $C, D\subset A_m$ with non abelian central projection such that
$$\|[\phi_{n,m}(f), c]\|<\varepsilon \|c\| ~~\mbox{and}~~ \|[\psi_{n,m}(f), d]\|<\varepsilon \|d\|~~~\mbox{for all} ~~~f\in F,~ c\in C,~ d\in D.$$
Consequently, both $C^*$algebras $A$ and $B$ are approximately divisible in the sense of Definition 1.2 of [BKR].  By Theorem 2.3 of [TW], both $A$ and $B$ are  ${\cal Z}$-stable. That is, $A\otimes {\cal Z}\cong A$ and $B\otimes {\cal Z}\cong B$, where ${\cal Z}$ is the Jiang-Su algebra (see [JS]). Furthermore, by using [Ti] (see [Cow-Ell-I] also), one can prove that $Cu(A)\cong Cu(B)$ and $Cu(A\otimes C(S^1))\cong Cu(B\otimes C(S^1))$.



$\\$ Guihua Gong, College of Mathematics and Information Science, Hebei Normal University,  Shijiazhuang, Hebei, 050024, China.\\
Department of Mathematics, University of Puerto Rico at Rio Piedras, PR 00936, USA\\
email address:  guihua.gong@upr.edu\\
Chunlan Jiang, College of Mathematics and Information Science, Hebei Normal University,  Shijiazhuang, Hebei, 050024, China.~~\\
email address:  cljiang@hebtu.edu.cn\\
Liangqing Li, Department of Mathematics, University of Puerto Rico at Rio Piedras, PR 00936, USA\\
email address:  liangqing.li@upr.edu\\

\clearpage

\begin{tiny}

\begin{small}

\end{small}

\end{tiny}

\end{document}